\documentclass[a4paper,11pt]{article}

\usepackage[mac]{inputenc}
 \usepackage[T1]{fontenc}
 \usepackage[normalem]{ulem}
 \usepackage[english]{babel}
\usepackage{amsmath}
  \usepackage{amsthm}
 \usepackage{bbm}
 \usepackage{amssymb}
 \usepackage{verbatim}
 \usepackage{vmargin}
 \usepackage{graphicx}
 \usepackage{color}
 \usepackage{epsfig}
 \usepackage[only,llbracket,rrbracket]{stmaryrd}
 \usepackage{enumerate}

\usepackage{setspace}

\newtheorem{thm}{Theorem}[section]

\newtheorem*{notat}{Notations}
\newtheorem{prop}{Proposition}[section]
\newtheorem{lem}{Lemma}[section]

\numberwithin{equation}{section}

\theoremstyle{theorem}
{\vskip 0.5cm}

\newtheorem{rque}{\textbf{Remark}}[section]{\vskip 0.5cm} 
{\vskip 0.5cm}

\title{Quasineutral limit of the Vlasov-Poisson system \\ with massless electrons}
\author{Daniel Han-Kwan\footnote{\'Ecole Normale Sup\'erieure, D\'epartement de Math\'ematiques et Applications,  45 rue d'Ulm 75230 Paris Cedex 05 France, email : {daniel.han-kwan@ens.fr}}}
\date{}
\begin{document}
 \maketitle

\begin{abstract}
In this paper, we study the quasineutral limit (in other words the limit when the Debye length tends to zero) of  Vlasov-Poisson like equations describing the behaviour of ions in a plasma. We consider massless electrons, with a charge density following a Maxwell-Boltzmann law. For cold ions, using the relative entropy method, we derive the classical Isothermal Euler or the (inviscid) Shallow Water systems from fluid mechanics. In a second time, we study the combined quasineutral and strong magnetic field regime for such plasmas.
\end{abstract}


\section{Introduction}

\subsection{Physical motivation}

For high enough kinetic temperatures exceeding the atomic ionization energy, atoms tend to decompose into electrons and ions (that is, negatively and positively charged particles):  a plasma is a physical or chemical system where such a ionization has occured.
Roughly speaking, we simply consider that plasmas are gases made of positive and negative charges. Unlike gases, plasmas are highly conductive. As a consequence, particles interact with each other by creating their own electromagnetic fields which can dramatically affect their behaviour.

The plasma state is considered as the fourth state of matter. Actually it is the most common one in the universe : it is widely recognized that at least $95\%$ of the matter consists of it ! For instance, the suns and other stars are filled with plasma, so is the interstellar medium and so on.  Terrestrial plasmas are also quite easy to find: they appear in flames, lightning or in the ionosphere. For the last decades, there has been an increasing interest in creating artificial plasmas, for experimental or industrial purposes. For instance, neon lamps or plasma displays for televisions are now part of our everyday life. An extremely promising application of plasmas consists in the fusion energy research (by magnetic or inertial confinement). This paper specifically aims at rigorously deriving some mathematical models which would help to understand the physics in tokamaks, which are the boxes in which plasmas from magnetic confinement fusion are contained.


 \subsubsection{Basic kinetic models for plasmas}

We adopt a statistical description of the plasma: we describe the behaviour of the charged particles by considering kinetic equations satisfied by their repartition function. That means that we do not follow particles one by one by solving Newton equations but are rather interested in their collective behaviour. 

We present the mathematical models we are going to study in the following. In order to establish them, we have to make some standard approximations which we now explain.

 \begin{itemize}
\item \textbf{Assumption 1}:  We assume the plasma to be collisionless. Thus, we will consider Vlasov-like equations without collision operators.
\item\textbf{Assumption 2}: The plasma is non-relativistic and the electric field $E$ is electrostatic. This means that we consider electromagnetic fields that satisfy the electrostatic approximation of the Maxwell equations:
\begin{equation*}
\left\{
    \begin{array}{ll}
    \operatorname{rot} E =0, \\
    \operatorname{div} E = \frac{\rho}{\epsilon_0}.
        \end{array}
  \right.
\end{equation*}
We can also consider an additional magnetic field $B$ satisfying
\[
 \operatorname{div} B = 0
\]
and it has to be stationary in time in order to be consistent with the electrostatic approximation. The electromagnetic fields act on the charged particles (with charge $q$) through the Lorentz force:
\[
F= q (E + v \wedge B).
\]

\item \textbf{Assumption 3}: The plasma evolves in a domain without boundaries. This means in particular that we may restrict to periodic data in the space variable, which may seem unrealistic, but which is commonly done in plasma physics and mathematics.

 \end{itemize}
 We discuss the physical relevancy of these assumptions in the following Remark:
 
 \begin{rque}
  \begin{itemize}
\item Discussion on \textbf{A1}: Let $\Lambda=4\pi N_0 \left(\frac{\epsilon_0 T }{N_0 e^2}\right)^{3/2}$, where $N_0$ is the average number density of particles, $\epsilon_0$  is the vacuum permittivity, $T$ the average temperature of the plasma and $e$ is the fundamental electric charge.

The typical collision frequency is given by:
$$\omega_c = \frac{e^4 \log \Lambda}{4\pi \epsilon_0^2m^{1/2}} \frac{N_0}{T^{3/2}}.$$ 
So the plasma can be considered as collisionless if it is diffuse and high temperature. Most plasmas can be considered as collisionless to a very good approximation \cite{Fitz}.
\item Discussion on \textbf{A2}:  The electrostatic approximation is relevant as soon as 
$$\frac{c\tau}{L}>>1.$$
denoting by $c$ the speed of light, $\tau$ the characteristic observation time and $L$ the characteristic observation length. Therefore this approximation often appears as reasonable in practical situations for terrestrial plasmas. At least, it is valid for short observation lengths.

 \item Discussion on \textbf{A3}: By ignoring boundary effects, we neglect some important physics, such as the formation of the Debye sheath near walls, which are boundary layers often surrounding plasmas confined in some material.  Very few seems to be mathematically known about this phenomenon, starting from a Vlasov-Poisson equation.

\end{itemize}
 \end{rque}
 
 Within these approximations, the kinetic system reads:
 
 \begin{equation}
 \label{VP}
\left\{
    \begin{array}{ll}
  \partial_t f_i + v.\nabla_x f_i+ q_i/ m_i \left(E+v\wedge B\right).\nabla_v f_i =0  \\
    \partial_t f_e + v.\nabla_x f_e + q_e /m_e \left(E+v\wedge B\right).\nabla_v f_e =0  \\
  E = -\nabla_x V\\
  -\Delta_x V = \frac{1}{\epsilon_0}\left(q_i \int f_i dv +q_e \int f_e dv\right),\\
    \end{array}
  \right.
\end{equation}
with $x \in \mathbb{R}^n \text{ or } \mathbb{T}^n = \left( \mathbb{R}/ \mathbb{Z} \right)^n, v \in \mathbb{R}^n, t \in \mathbb{R}^+$. 
We can associate to these equations the initial conditions:
 \begin{equation}
\left\{
    \begin{array}{ll}
    f_{i,\vert t=0} =f_{i,0} ,\quad f_{i,0}\geq 0, \quad  \int f_{i,0} dvdx =1,\\
       f_{e,\vert t=0} =f_{e,0},\quad f_{e,0}\geq 0, \quad  \int f_{e,0} dvdx =1.\\
    \end{array}
  \right.
\end{equation}

The parameter $n$ is the space dimension, equal to $1,2$ or $3$ in the following. Quantity $f_i$ (resp. $f_e$) is interpreted as the density distribution of ions (resp. electrons) : $f(t,x,v)dxdv$ is interpreted as the probability of finding particles at time $t$ with position $x$ and velocity $v$. 
The parameter $m_i$ (resp. $m_e$) is the mass of one ion (resp. electron). Likewise, $q_i$ (resp. $q_e$) is the charge of one ion (resp. electron). For simplicity we will take $q_e= -e$ and $q_i=e$ .

We now intend to reduce the two transport equations into only one.
To this end, we can observe that the mass ratio between electrons and ions is very small:
 $$\frac{m_e}{m_i} << 1,$$ 
so that qualitatively the two types of particles have really different dynamical behaviour. 
Therefore we make the additional approximation for our idealized model:
\begin{itemize}
\item \textbf{Assumption 4}: The mass ratio between ions and electrons is infinite:  $\frac{m_i}{m_e}= +\infty$.
\end{itemize}

This remark allows to reduce System (\ref{VP}) to only one transport equation. Depending on the interpretation of Assumption 4, we get two classes of models:

\subsubsection*{Infinite mass ions ($m_i=+\infty$)}

 One can consider the point of view of electrons, from which ions are very slow, motionless at  equilibrium:

\begin{equation}
n_i= \int f_i dv = 1.
\end{equation}

Then, assuming there is no magnetic field, system (\ref{VP}) written in dimensionless variables reduces to the classical Vlasov-Poisson system:
\begin{equation}
\left\{
    \begin{array}{ll}
  \partial_t f + v.\nabla_x f + E.\nabla_v f =0  \\
  E = -\nabla_x V\\
  -\Delta_x V = \int f dv -1\\
    f_{\vert t=0} =f_0\geq 0, \int f_0 dvdx =1.\\
    \end{array}
  \right.
\end{equation}

This system was intensively studied in the mathematical literature (for $n= 3$ in particular) and the Cauchy-Problem is rather well understood. We refer to the works of Arsenev \cite{Ar}, Horst and Hunze \cite{HH} for global weak solutions, DiPerna and Lions \cite{DPL} for global renormalized solutions, Pfaffelmoser \cite{Pfa}, Schaeffer \cite{Scha} for classical solutions,  Lions and Perthame \cite{LP} for weak solutions with high order velocity moments and  Loeper \cite{Loe} on the uniqueness problem.

\subsubsection*{ Massless electrons ($m_e=0$)}

 Otherwise, one can consider the viewpoint of ions : electrons then move very fast and quasi-instantaneously reach their local thermodynamic equilibrium\footnote{Actually since $m_e<<m_i$, the typical collision frequency for the electrons is much larger than for the ions and thus collisions for the electrons may be not negligeable: for this reason they can reach their local thermodynamic equilibrium.}. Then their density $n_e$ follows the classical Maxwell-Boltzmann law (see \cite{LL}) :
\begin{equation}
\label{MB}
n_e = \int f_e dv = d(x)\exp \left(\frac{eV}{k_B T_e}\right),
\end{equation}
 where $V$ denotes the electric potential, $k_B$ is the Boltzmann constant, $T_e$ the average temperature of the electrons,  $d \in L^1(\mathbb{R}^n)$ is a term due to an external potential preventing the particles from going to infinity (we also refer to \cite{Bou} and references therein). 
 
More precisely, we have:
\begin{equation}
d(x)=n_0 e^{\frac{-H(x)}{k_B T_e}}, 
\end{equation}
where $n_0 \in \mathbb{R}$ is a normalizing constant and $H$ is the external confining potential.
 
 The Poisson equation then reads:
 \begin{equation}
 -\Delta_x V = \int f dv -{d\exp \left(\frac{eV}{k_B T_e}\right)
}.
\end{equation}

 One should notice that in this case, in general
 $$\int \left(\int f dv -d\exp \left(\frac{eV}{k_B T_e}\right)\right) dx \neq 0, $$
meaning that  global neutrality does not hold, since the total charge of electrons is not a priori fixed\footnote{This feature will prevent us from studying this system on the torus; instead we will do so on the whole space $\mathbb{R}^3$.}. 

We may also consider the case when the total charge of the electrons is fixed, in which case the Poisson equation reads:
\begin{equation}
 -\Delta_x V = \int f dv -\frac{d\exp \left(\frac{eV}{k_B T_e}\right)
}{\int_{\mathbb{R}^n} d \exp \left(\frac{eV}{k_B T_e}\right)
 dx}.
\end{equation}

The existence of global weak solutions to these two systems in dimension three has been investigated by Bouchut \cite{Bou}. We will recall some of the properties of these solutions in Section \ref{GW}.

  An approximation widely used in plasma physics consists in linearizing the exponential law:
 \begin{equation}
 n_e=d\left(1+\frac{eV}{k_B T_e}\right).
 \end{equation}
 This approximation is valid  from the physical point of view as long as:
 $$\frac{eV}{k_B T_e} <<1,$$
 that is as long as the electric energy is small compared to the kinetic energy.
 
 We will consider this law in the case of the torus $\mathbb{T}^n$ (with $n=1,2$ or $3$), thus we do not need a confining potential (and we take $d=1$).
 
 In the following we will only focus on models with such Maxwell-Boltzmann laws.

\subsubsection{The Debye length}
\label{debye}
We define now the Debye length $\lambda_D^{(\alpha)}$ as:

\begin{equation}
\lambda_D^{(\alpha)} = \sqrt{\frac{\epsilon_0 k_B T_\alpha}{N_\alpha e^2}},
\end{equation}
where $k_B$ is the (universal) Boltzmann constant, $T_\alpha$ and $N_\alpha$ are respectively the average temperature and density of electrons (for $\alpha=e$) or ions (for $\alpha=i$).

 The Debye length is a fundamental parameter which is of tremendous importance in plasmas. It can be interpreted as the typical length below which charge separation occurs. 

In plasmas, this length may vary by many orders of magnitude (Typical values go from $10^{-3}m$ to  $10^{-8}m$). In practical situations, for terrestrial plasmas, it is always small compared to the other characteristic lengths under consideration, in particular the characteristic observation length, denoted by $L$. Actually, the condition $\lambda_D <<L$ is sometimes required in the definition itself of a plasma.

Therefore, if we set:
$$\frac {\lambda_D} {L} = \epsilon<<1,$$
then in many regimes, it is relevant, after considering relevant dimensionless variables, to consider that the Poisson equation formally reads :
$$
-\epsilon^2 \Delta_x V_\epsilon = \pm \left(n_i - n_e\right).
$$

The quasineutral limit precisely consists in considering the limit $\epsilon\rightarrow 0$.

\subsubsection{Why quasineutral fluid limits ?}

From the numerical point of view, kinetic equations are harder to handle than fluid equations. Indeed the main difficulty is that we have to deal with a phase space of dimension $6$ (for $x,v\in \mathbb{R}^3$). Actually, another outstanding problem for simulating plasmas is the following : there are characteristic lengths and times of completely different magnitude (think of the Debye length and the observation length) that make numerics really delicate.

In this work, we particularly aim at getting simplified hydrodynamic systems after taking quasineutral limits.
Simplified fluid models have some advantages:

\begin{itemize}

\item With a fluid description, we deal with a phase space of lesser dimension. Furthermore after taking the limit we now handle only one characteristic time and length. For these reasons, numerical simulations are easier to perform. 

Of course it is well-known that the fluid approximation is not always accurate for simulations of plasmas, but it is nevertheless valid in some regimes that we may describe in the analysis. So it is important to be aware of the physical assumptions we make when we derive the equations.

\item Macroscopic quantities, such as charge density or current density are easier to experimentally measure (by opposition, the repartition function is out of reach). So this is a way to check if the initial modeling is accurate or not.

\item A simplified fluid description can help us to qualitatively describe the behaviour of the plasma.
\end{itemize}

The derivation of limit models is deeply linked to issues related to the research project of magnetic confinement fusion. For the last few years, there has been a wider interest in finding simplified systems to model quasineutral plasmas for devices such as tokamaks. Therefore a good mathematical understanding of these becomes important, as it would establish some theoretical basis to compare various models, like gyrokinetic, gyrofluid, MHD or Euler-like equations and understand their range of validity.

In this paper, we will focus on two kinds of quasineutral problems. In both problems, the starting point is the Vlasov-Poisson system with massless electrons (in other words with electrons following a Maxwell-Boltzmann law). First we will investigate the quasineutral limit alone, then we will in addition consider a large magnetic field and study the behaviour of the plasma in this regime.

\subsection{Quasineutral limit of the Vlasov-Poisson system with massless electrons}

In the first place, we are interested in the quasineutral limit for Vlasov-Poisson systems with Boltzmann-Maxwell laws and without magnetic field.
We will focus in particular on the limit $\epsilon \rightarrow 0$ for the following dimentionless system (we refer to the Annex for details on the scaling):

 \textbf {System (S)}  : Maxwell-Boltzmann law
 (for $x \in \mathbb{R}^3 \text, v \in \mathbb{R}^3, t \in \mathbb{R}^+$)
 
  \begin{equation}
  \label{(S)}
\left\{
    \begin{array}{ll}
  \partial_t f_\epsilon + v.\nabla_x f_\epsilon + E_\epsilon.\nabla_v f_\epsilon =0  \\
  E_\epsilon = -\nabla_x V_\epsilon\\
-\epsilon \Delta_x V_\epsilon = \int f_\epsilon dv - d e^{V_\epsilon}\\
    f_{\epsilon,\vert t=0} =f_{0,\epsilon} \geq 0, \int f_{0,\epsilon} dvdx =1.\\
    \end{array}
  \right.
\end{equation}

The method of proof we follow allows also to treat the case of variants of system (S), so we will mention the results we can get, without providing complete proofs, for the systems:

\begin{itemize}

\item \textbf{System (S')} : Maxwell-Boltzmann law with fixed total charge
 (for $x \in \mathbb{R}^3, v \in \mathbb{R}^3, t \in \mathbb{R}^+$)
  \begin{equation}
    \label{(S')}
\left\{
    \begin{array}{ll}
  \partial_t f_\epsilon + v.\nabla_x f_\epsilon + E_\epsilon.\nabla_v f_\epsilon =0  \\
  E_\epsilon = -\nabla_x V_\epsilon\\
-\epsilon \Delta_x V_\epsilon = \int f_\epsilon dv - \frac{d e^{V_\epsilon}}{\int d e^{V_\epsilon}}\\
    f_{\epsilon,\vert t=0} =f_{0,\epsilon} \geq 0, \int f_{0,\epsilon} dvdx =1.\\
    \end{array}
  \right.
 \end{equation}
  
\item \textbf {System (L)} :  Linearized Maxwell-Boltzmann law
 (for $x \in \mathbb{T}^n , v \in \mathbb{R}^n, t \in \mathbb{R}^+$ and $n=1,2,3$)

 \begin{equation}
   \label{(L)}
\left\{
    \begin{array}{ll}
  \partial_t f_\epsilon + v.\nabla_x f_\epsilon + E_\epsilon.\nabla_v f_\epsilon =0  \\
  E_\epsilon = -\nabla_x V_\epsilon\\
 V_\epsilon -\epsilon \Delta_x V_\epsilon = \int f_\epsilon dv -1\\
    f_{\epsilon,\vert t=0} =f_{0,\epsilon} \geq 0, \int f_{0,\epsilon} dvdx =1.\\
    \end{array}
  \right.
\end{equation}

\end{itemize}

For systems (S) and (S') we will from now on assume the boundedness properties on $d(x)=e^{-H(x)}$:
\begin{eqnarray}
\label{assuH-1}d=e^{-H} \in L^1\cap L^\infty(\mathbb{R}^3).\\
\label{assuH-2}\nabla_x H \in W^{s,\infty} \quad \text{for any} \quad s\in \mathbb{N}.
\end{eqnarray}
For instance, this holds for $H(x)= \sqrt{1+ \vert x \vert^2}$.

\begin{rque} From the mathematical viewpoint, we have to add the confining potential $H$ to ensure that the local density of electrons belongs to $L^1(\mathbb{R}^3)$.
\end{rque}

\subsubsection{Formal derivation of the isothermal Euler system from systems (S) and (S')}

 We will prove the local in time strong convergence of the charge density and current density:
$$\left( \rho_\epsilon:= \int f_\epsilon dv, J_\epsilon:=\int f_\epsilon v dv \right)$$  
to the local strong solution $(\rho,\rho u)$ to some Euler-type system, for initial data close (in some sense to be made precise later) to monokinetic data that is, 
$$f_\epsilon(t,x,v) \sim \rho_\epsilon(t,x) \delta_{v=u_\epsilon(t,x)},$$
with $u_\epsilon =\frac{J_\epsilon}{\rho_\epsilon}$.

Let us show now how we can guess what is the limit system. First, by integrating the Vlasov equation against $1$ and $v$, we straightforwardly get the local conservation laws  satisfied by the first two moments.
\begin{eqnarray}
\partial_t \rho_\epsilon + \nabla_x. J_\epsilon=0, \\
\partial_t J_\epsilon +  \nabla_x : \left(\int v \otimes v f_\epsilon dv \right) = \rho_\epsilon E_\epsilon.
\end{eqnarray}

Let us directly consider monokinetic data, i.e. $f_\epsilon(t,x,v) = \rho_\epsilon (t,x) \delta_{v=u_\epsilon(t,x)}$. The local conservations laws reduce to:
\begin{equation}
\partial_t \rho_\epsilon + \nabla_x. (\rho_\epsilon u_\epsilon)=0.
\end{equation}

\begin{equation}
\partial_t (\rho_\epsilon u_\epsilon) +  \nabla_x : (\rho_\epsilon u_\epsilon \otimes u_\epsilon) = -\rho_\epsilon \nabla_x V_\epsilon.
\end{equation}

In the case of (S) the Poisson equation reads:
\begin{equation*}
 -\epsilon \Delta_x V_\epsilon = \int f_\epsilon dv - de^{V_\epsilon}.
\end{equation*}

Since $\rho_\epsilon$ and $J_\epsilon$ are uniformly bounded in $L^\infty_t(L^1_x)$, the following convergences hold (up to a subsequence) in the sense of distributions: $\rho_\epsilon \rightharpoonup \rho$ and $J_\epsilon \rightharpoonup J$.

If we formally pass to the limit $\epsilon \rightarrow 0$ we get:
\begin{equation}
de^V=\rho.
\end{equation}

Consequently, we have $V= \log \left(\rho/d\right) $  and therefore $-\nabla_x V = - \frac{\nabla_x \rho}{\rho} + \frac{\nabla_x d}{d}$. Notice that $  \frac{\nabla_x d}{d} = -\nabla_x H$. Thus, the asymptotic equation we can expect is the following compressible Euler-type model (which can be interpreted as the isothermal Euler equation with an external confining force):

\begin{equation}
\label{euler}
\left\{
\begin{array}{ll}
\partial_t \rho + \nabla_x. (\rho u)=0, \\
\partial_t u +  u.\nabla_x u = - \frac{\nabla_x \rho}{\rho} -  \nabla_x H.

\end{array}
\right.
\end{equation}

In the case of (S') the Poisson equation reads:
\begin{equation}
-\epsilon \Delta_x V_\epsilon = \int f_\epsilon dv - \frac{d e^{V_\epsilon}}{\int d e^{V_\epsilon}dx}.
\end{equation}

If we formally pass to the limit $\epsilon \rightarrow 0$ we get:
\begin{equation}
 \frac{d e^{V}}{\int d e^{V}dx}= \rho.
\end{equation}

Consequently we have $V= \log \left(\frac{ \rho}{d} \int d e^{V}dx\right)$ and so, we get the same Euler equation (\ref{euler}).

\begin{rque}(Physical signification of monokinetic data)

We define:
$$T_{i,\epsilon}=\frac{1}{3\rho_\epsilon}\int f _\epsilon \left\vert v- \frac{J_\epsilon}{\rho_\epsilon} \right\vert^2 dv dx.$$ 

The quantity $T_{i,\epsilon}$ is nothing but the scaled temperature of the ions. 

Considering monokinetic data corresponds to the "cold ions" assumption, that is:
$$\int f \left\vert v- \frac{J}{\rho} \right\vert dv dx =0,$$ 
which means that we consider that the temperature of ions is equal to $0$. 

More precisely, the cold ions approximation means from the physical point of view that 
$$T_i << T_e.$$ 
It turns out that this approximation is highly relevant for terrestrial plasmas  and widely used in plasma physics, especially in the study of tokamak plasmas.

There are two main physical reasons why it is relevant to consider that the temperature of electrons is much higher than the temperature of ions : first of all , there exist many plasma sources which can heat the electrons more strongly than the ions. Second notice that energy transfer in a two-body collision is much more efficient if the masses are similar. Thus, since ions and electrons have very different masses, there is almost no transfer of energy from the electrons to the ions. For instance this approximation is used in order to derive the classical Hasegawa-Mima equation (\cite{HM1}).

\end{rque}
\begin{rque}Nevertheless we observe in the isothermal Euler limit system (\ref{euler}) that the ions evolve as if they had the temperature of electrons (of order $1$) ! Moreover, ions seem to have better confinement properties than expected, since they feel the confining potential in the limit equation.

\end{rque}

\begin{rque}  For $(L)$ the corresponding Euler-type system is the following:
\begin{equation}
\label{shallow}
\left\{
\begin{array}{ll}
\partial_t \rho + \nabla_x. (\rho u)=0, \\
\partial_t u +  u.\nabla_x u = - \nabla_x \rho.
\end{array}
\right.
\end{equation}

Actually System (\ref{shallow}) can be interpreted in $1D$ or $2D$ as an inviscid Shallow Water system (and $\rho$ is then understood as the depth of the fluid). This quite remarkable fact is one amongst many analogies between geophysics and plasma physics models (see for instance the work of Hasegawa and Mima \cite{HM1} and the review paper \cite{DA}). For instance the concept of "zonal flows" is used in both fields and the mechanism responsible for their generation may be the same. Only the name differs: drift waves for plasma physics, Rossby waves for geophysics, see Cheverry, Gallagher, Paul and Saint-Raymond for a recent mathematical study \cite{CGPSR}.
\end{rque}

\subsubsection{Principle of  the proof : the relative entropy method}

The relative entropy method (also referred to as the modulated energy method) was first introduced in kinetic theory independently by Golse \cite{BGP} in order to study the convergence of solutions to a scaled Boltzmann equation to solutions of incompressible Euler for well-prepared data and some technical assumptions (see Saint-Raymond \cite{SR} for the latest developments on the topic) and by Brenier \cite{Br} in order to derive incompressible Euler equations from the quasineutral Vlasov-Poisson equation for electrons in a fixed background of ions.

More precisely Brenier shows the convergence as $\epsilon \rightarrow 0$ of the first two moments $(\rho_\epsilon, J_\epsilon):=(\int f_\epsilon dv, \int v f_\epsilon dv)$ of the starting system:

\begin{equation}
\left\{
    \begin{array}{ll}
  \partial_t f_\epsilon + v.\nabla_x f_\epsilon + E_\epsilon.\nabla_v f_\epsilon =0  \\
  E_\epsilon = -\nabla_x V_\epsilon\\
  -\epsilon\Delta_x V_\epsilon = \int f_\epsilon dv -1\\
    f_{\epsilon, \vert t=0} =f_{0,\epsilon} \geq 0, \int f_{0,\epsilon} dvdx =1.\\
    \end{array}
  \right.
\end{equation}
to the smooth solution of the limit system which is the classical incompressible Euler system:
\begin{equation}
\left\{
    \begin{array}{ll}
    \rho=1 \\
  \partial_t u  + u.\nabla_x u + \nabla_x p=0  \\
\operatorname{div} u=0.\\
    \end{array}
  \right.
\end{equation}

Brenier treated the case of well-prepared monokinetic data (i.e. cold electrons); the convergence was then generalized by Masmoudi \cite{Mas} for ill-prepared monokinetic data. 

We mention the works  \cite{BOS}, \cite{BL}, \cite{BMP}, \cite{GSR2}, \cite{GJV}, \cite{PSR} which also use the relative entropy method in order to derive fluid equations from Vlasov-like systems.

The principle of the method is the following. For system (S), it can be shown that the following functional is non-increasing:
\begin{equation}
\mathcal{F}_\epsilon(t)= \frac{1}{2}\int f_\epsilon\vert v \vert ^2 dv dx + \frac{\epsilon}{2} \int \vert \nabla_x V_\epsilon \vert^2 dx
 + \int d(x)(V_\epsilon-1)e^{V_\epsilon} dx.  \end{equation}
We call this functional the energy of the system.

We then consider the functional $\mathcal{H}_\epsilon$ which is built as a modulation of this energy:
\begin{equation}
\mathcal{H}_\epsilon(t)=\frac{1}{2}\int f_\epsilon\vert v - u\vert ^2 dv dx + \frac{\epsilon}{2} \int \vert \nabla_x V_\epsilon \vert^2 dx+ \int (de^{V_\epsilon} \log \left( de^{V_\epsilon} /\rho \right)  - de^{V_\epsilon} +\rho)dx ,
\end{equation}
with $(\rho,u)$ a smooth solution to Isothermal Euler.

Quite surprisingly, it turns out that the last term of $\mathcal{H}_\epsilon$ is similar to the usual relative entropy for collisional (such as Boltzmann or BGK) equations.

What we want to prove is that this functional is in fact a Lyapunov functional. We will show that indeed, $\mathcal{H}_\epsilon$ satisfies some stability estimate:
\begin{eqnarray}
{\mathcal{H}}_\epsilon(t) \leq   {\mathcal{H}}_\epsilon(0) +  G_\epsilon(t) +  C\int_0^t \Vert \nabla_x u \Vert_{L^\infty} \mathcal{H}_\epsilon(s) ds,
\end{eqnarray}
with $G_\epsilon(t) \rightarrow_{\epsilon\rightarrow 0} 0$ uniformly in time.

Then, assuming that the initial conditions is well-prepared in the sense that
$$\mathcal{H}_\epsilon(0) \rightarrow_{\epsilon\rightarrow 0} 0,$$
this yields that $\mathcal{H}_\epsilon(t) \rightarrow_{\epsilon\rightarrow 0} 0$, thus proving the strong convergence in some sense (which will be made precise later on) to smooth solutions of the isothermal Euler equation, as long as the latter exist. The proof relies on the fine algebraic structure  of  the nonlinearities in systems (S) and (S').
One major advantage of this method is that it only requires weak regularity on the solutions to the initial system but allows to prove limits in a strong sense. Nevertheless, it requires a good understanding of the Cauchy problem for the limit system (in particular, we must have a notion of stability for the limit). It should be noticed that even if we considered very smooth solutions  (say for instance $H^s$ with $s$ large) to the initial system, we would not be able to propagate uniform bounds and thus prove compactness. Indeed, the only uniform controls we have are the energy bound and the conservation of $L^p$ norms of the number density.

Basically this is nothing but a stability result : roughly speaking , this result tells us that monokinetic solutions are stable with respect to perturbations of the energy.

Let us mention that the method used for these two systems can also apply to the quasineutral limit for an isothermal Euler-Poisson version of system (S) studied by Cordier and Grenier \cite{CG}. We refer to Section \ref{IEP}.

\subsection{Quasineutral limit for the Vlasov-Poisson equation with massless equations and with a strong magnetic field}

Next we are interested in the behaviour of the plasma if one applies an intense magnetic field. Such a regime is particularly relevant for plasmas encountered in magnetic confinement fusion research. Plasmas are expected to be confined inside tokamaks thanks to this magnetic field. One challenging mathematical problem is to rigorously prove if this strategy is likely to succeed or not.

In this paper, we consider the simplest geometric case of a constant magnetic field with a fixed direction and a fixed (large) intensity.

\subsubsection{Scaling of the Vlasov equation}

We first introduce some notations:
\begin{notat}

\begin{itemize}
Let $(e_1,e_2,e_\parallel)$ be a fixed orthonormal basis of $\mathbb{R}^3$.
\item The subscript $\perp$ stands for the orthogonal projection on the plane $(e_1,e_2)$, while the subscript $\parallel$ stands for the projection on $e_\parallel$ .
\item For any vector $X=(X_1,X_2,X_\parallel)$, we define $X^\perp$ as the vector $(X_y,-X_x,0)=X\wedge e_\parallel$.
\end{itemize}
\end{notat}

We consider a strong magnetic field of the form:
$$B = B e_\parallel.$$
Roughly speaking, ``strong'' means that $\vert B \vert \sim 1/\epsilon$.
With this time a quasineutral ordering of the form
\[
\frac{\lambda_D}{L} = \epsilon^\alpha,
\]
where $\alpha>0$ is an arbitrary parameter, we get in the end the quasineutral system:

\begin{equation}
\label{(S')-B}
\left\{
    \begin{array}{ll}
  \partial_t f_{\epsilon} + v.\nabla_x f_{\epsilon} + \left(E_{\epsilon}+ \frac{v\wedge e_\parallel}{\epsilon}\right).\nabla_v f_{\epsilon} =0  \\
  E_{\epsilon} = -\nabla_x V_{\epsilon}\\
-\epsilon^{2\alpha} \Delta_x V_{\epsilon} = \int f_{\epsilon} dv -\frac{d e^{V_{\epsilon}}}{\int d e^{V_\epsilon}dx}
\\
    f_{{\epsilon},\vert t=0} =f_{0,\epsilon}, \quad \int f_{0,\epsilon} dvdx =1.\\
    \end{array}
  \right.
\end{equation}
We refer to the Annex for details on the scaling. The range of parameters $\alpha>1$ is particularly relevant from the physical point of view.

\subsubsection{Comments on the expected result}

We will study the limit, once again by using the relative entropy method and show the convergence of the first two moments $(\rho_\epsilon, u_\epsilon)$ (defined as before) to smooth solutions to the system:

\begin{equation}
\left\{
    \begin{array}{ll}
\partial_t \rho + \partial_{x_\parallel} (\rho w_\parallel) =0, \\
\partial_t w + w_\parallel \partial_{x_\parallel} w = -\frac{\nabla_{x_\parallel} \rho}{\rho} - \nabla_{x_\parallel} H.
    \end{array}
  \right.
\end{equation}

We observe that there is no more dynamics in the orthogonal plane (that is, in the $x_\perp$ variable), which can be interpreted as a good confinement result.


For the study of this limit we will have to face more technical difficulties than without magnetic field. Indeed, the strong magnetic field engenders time oscillations of order $\mathcal{O}(1/\epsilon)$ on the number density. Consequently in order to show strong convergence we will have to:
\begin{itemize}
\item   filter out the time oscillations.
\item  add some correction of order $\mathcal{O}(\epsilon)$ to the limit $(\rho,w)$ in order to get an approximate zero of the so-called acceleration operator (\ref{acc-op}).
\end{itemize}

The striking point here is that we can study the limit for any value of $\alpha$. In contrast, for the system describing the electrons with heavy ions:

$$
\left\{
    \begin{array}{ll}
  \partial_t f_{\epsilon} + v.\nabla_x f_{\epsilon} + \left(E_{\epsilon}+ \frac{v\wedge e_\parallel}{\epsilon}\right).\nabla_v f_{\epsilon} =0  \\
  E_{\epsilon} = -\nabla_x V_{\epsilon}\\
-\epsilon^{2\alpha}\Delta_x V_{\epsilon} = \int f_{\epsilon} dv -1\\
    f_{{\epsilon},\vert t=0} =f_{0,\epsilon}, \quad \int f_{0,\epsilon} dvdx =1.\\
    \end{array}
  \right.
$$
it seems primordial to take $\alpha = 1$, so that the Debye length and the Larmor radius vanish at the same rate. This specific scaled system was studied by Golse and Saint-Raymond in \cite{GSR2}.

The heuristic underlying reason is that the Poisson equation with a Maxwell-Boltzmann law is in some sense more stable in the quasineutral limit than the "usual" one. Indeed the electric potential is in the limit explicitly a function of $\rho$, whereas in the "usual" case, it appears as a Lagrange multiplier or equivalently as a pressure.


\begin{rque}
We could as well set $\frac{\lambda_D}{L}= \delta$ and let $\epsilon,\delta$ go to $0$ independently. One can readily check that we would get the same results.
\end{rque}

\begin{rque}
We may also consider the linearized Maxwell-Boltzmann law and perform the same convergence analysis.
\end{rque}

\subsection{Outline of the paper}

The following of this article is organized as follows. First in section \ref{GW}, we recall some elements on the global weak solutions theory for systems (S), (S') and (L) and recall some useful \textit{a priori} uniform bounds.
This theory is due to Arsenev \cite{Ar} and Bouchut \cite{Bou}. Then, in the spirit of Lions and Perthame \cite{LP}, we consider global weak solutions satisfying some convenient local conservation laws.

In section \ref{quasi-S}, we will focus on the quasineutral limit from the Vlasov-Poisson System (S) to an isothermal Euler system, using the relative entropy method (Theorem \ref{theoS}). 
The crucial step is to show the algebraic identity (\ref{crucial}) that describes the decay of the relative entropy and from which we will be able to get a stability inequality.

In section \ref{quasi-others}, we generalize the method for systems (S') and (L) (Theorems  \ref{theoS'} and \ref{theoL}), by only sketching the proofs. We will show that this method can also be applied to the quasineutral limit of a system previously studied by Cordier and Grenier \cite{CG}.

 Finally, in section \ref{quasi-magn}, we investigate the combined quasineutral and large magnetic field regime for system (S). The convergence result is stated in Theorem \ref{theoS'rot}.



\section{Global weak solutions and local conservation laws for the Vlasov-Poisson systems}
\label{GW}

\subsection{Global weak solutions theory}

Following Arsenev \cite{Ar}, it is straightforward to build global weak solutions to System (L), for any fixed $\epsilon>0$.
\begin{thm}
Let $n=1,2$ or $3$. We consider the functional:
\begin{equation}
\mathcal{E}_\epsilon(t)= \frac{1}{2}\int f_\epsilon\vert v \vert ^2 dv dx + \frac{1}{2}\int V_\epsilon^2 dx  + \frac{\epsilon}{2} \int \vert \nabla_x V_\epsilon \vert^2 dx. 
\end{equation}
For any $\epsilon>0$ and initial data $f_{0,\epsilon}\geq 0 $ bounded in $L^1 \cap L^\infty(\mathbb{R}^{2d}) $ such that $\mathcal{E}_\epsilon(0)$ is finite, there exists a global weak solution to (L) with $f_\epsilon \in L^\infty_t(L^1_{x,v}) \cap L^\infty_{t,x,v}$ and $\mathcal{E}_\epsilon(t)$ is non-increasing.
\end{thm}

Following Bouchut \cite{Bou}, we obtain the existence of global weak solutions to system (S) and (S').
We recall that $d$ satisfies assumptions (\ref{assuH-1}-\ref{assuH-2}) (in particular, $d\in L^1(\mathbb{R}^3)$).

\begin{thm}\label{exiS} For any $\epsilon>0$ and initial data $f_{0,\epsilon}\geq 0 $ bounded in $L^1 \cap L^\infty(\mathbb{R}^6) $ and satisfying  $\int (1 +\vert x \vert^2 + \vert v \vert^2) f_{0,\epsilon} dxdv<\infty$:
\begin{itemize}
\item \textbf{The case of (S)}
Let  $\mathcal{F}_\epsilon(t)$ be the functional defined as follows:
\begin{equation}
\mathcal{F}_\epsilon(t)= \frac{1}{2}\int f_\epsilon\vert v \vert ^2 dv dx + \int d(x)(V_\epsilon-1)e^{V_\epsilon} dx  + \frac{\epsilon}{2} \int \vert \nabla_x V_\epsilon \vert^2 dx.
\end{equation}
If $\mathcal{F}_\epsilon(0)$ is finite, there exists $f_\epsilon \in L^\infty_t(L^1_{x,v}) \cap L^\infty_{t,x,v}$ global weak solution to (S) with  $\mathcal{F}_\epsilon(t)$ non-increasing.

\item \textbf{The case of (S')}:
Let  $\mathcal{G}_\epsilon(t)$ be the functional defined as follows:
\begin{equation}
\begin{split}
\mathcal{G}_\epsilon(t)=& \frac{1}{2}\int f_\epsilon\vert v \vert ^2 dv dx + \int d(x)\left(V_\epsilon - \log \left(\int d e^{V_\epsilon} dx\right)\right)\frac{e^{V_\epsilon}}{\int d e^{V_\epsilon} dx} dx \\
   +& \frac{\epsilon}{2} \int \vert \nabla_x V_\epsilon \vert^2 dx.
   \end{split}
\end{equation}

If $\mathcal{G}_\epsilon(0)$ is finite,there exists $f_\epsilon \in L^\infty_t(L^1_{x,v}) \cap L^\infty_{t,x,v}$ global weak solution to (S') with  $\mathcal{G}_\epsilon(t)$ non-increasing.

\end{itemize}

In addition, in both cases, we have : $V_\epsilon \in L^\infty_t(L^6_x)$ and $\operatorname{ess sup}_{t,x} V_\epsilon <\infty$. In particular it means that $e^{V_\epsilon} \in L^\infty_{t,x}$.

\end{thm}


The main difficulties in \cite{Bou} are to get estimates for the electric potential in the Marcinkiewicz space $M^3$  to provide some strong compactness, and to use a relevant regularization scheme to preserve the energy inequality.

We assume from now on that the initial data satisfy the uniform estimates :
\begin{eqnarray} 
\label{assu1} \forall \epsilon>0, f_{0,\epsilon}\geq 0, \\
\forall \epsilon>0, f_{0,\epsilon} \in L^1\cap L^\infty(\mathbb{R}^n),  \text{   uniformly in   } \epsilon, \\
\exists C>0, \forall \epsilon>0, \quad \mathcal{E}_\epsilon(0)\leq C \quad (\text{resp.} \quad  \mathcal{F}_\epsilon, \quad \mathcal{G}_\epsilon).
\end{eqnarray}

Using a very classical property for Vlasov equations with zero-divergence in $v$ force fields, we get the following unifom in $\epsilon$ estimates.
\begin{lem} For $f_\epsilon$ global weak solution of (L) (resp. (S), resp. (S') we have
\begin{itemize}
\item (Conservation of $L^p$ norms) For any $p \in [1, +\infty]$, for any $t\geq 0$, $\Vert f_\epsilon(t) \Vert_{L^p_{x,v}} \leq \Vert f_\epsilon(0) \Vert_{L^p_{x,v}}$.
\item (Maximum principle) If $f_\epsilon(0) \geq 0$ then for any $t\geq 0, f_\epsilon (t)\geq 0$.
\item (Bound on the energy) $\forall t\geq 0$, $\mathcal{E}_\epsilon(t)\leq C$ (resp. $\mathcal{F}_\epsilon$, resp. $\mathcal{G}_\epsilon$).
\end{itemize}
\end{lem}

\begin{lem}Define $J_\epsilon(t,x)= \int f_\epsilon v dv$. Then $J_\epsilon \in L^\infty_t(L^1_x)$ uniformly with respect to $\epsilon$.
\end{lem}

\begin{proof}
Actually by the same method, we can also prove that $J_\epsilon  \in L^\infty_t(L^p_x)$ for some $p>1$ depending on the space dimension, but this result is sufficient for our purpose. The proof is very classical. We can first notice that there exists $C>0$ independent of $\epsilon$, such that:
$$
\int f_\epsilon \vert v \vert ^2 dv dx \leq C.
$$
For (L) this is clear by conservation of the energy since all the terms are non-negative.

In the case of (S) we observe
\begin{eqnarray*}
\int f_\epsilon \vert v \vert^2 dv dx &\leq& \mathcal{F}_\epsilon(t) - 2 \int d(V_\epsilon -1 ) e^{V_\epsilon} dx\\
&\leq&  \mathcal{F}_\epsilon(0) +2 \Vert d \Vert_{L^1},
\end{eqnarray*}
since for any $x \in \mathbb{R}, (x-1)e^x \geq - 1$. The case of (S') is of course similar.
Then we can simply write by positivity of $f_\epsilon$:
$$
\left\vert \int f_\epsilon v dv \right\vert \leq \int f_\epsilon \vert v \vert  dv \leq \int_{\vert v \vert \leq 1} f_\epsilon dv + \int_{\vert v \vert \geq 1} f_\epsilon \vert v\vert^2dv,
$$
so that: $\Vert J_\epsilon\Vert_{L^\infty_t(L^1_x)} \leq 1+ C$.

\end{proof}

\subsection{Local conservation laws}

We now assume the following additional non-uniform estimate on the initial data:
\begin{eqnarray}
\label{assu4}\exists m>3, \forall \epsilon >0, \int \vert v \vert^m f_{0,\epsilon}dvdx <\infty. 
\end{eqnarray}


This  allows to build  global weak solutions with high order moments in the spirit of Lions and Perthame \cite{LP}. We are then able to deal with solutions satisfying  local conservations  laws for charge and current.

For system (L), whose Poisson equation only differs from the usual one by the linear term $V$, we can straightforwardly use the results of \cite{LP}. We obtain global weak solutions with high order moments in $v$, known to satisfy as well the following local conservation laws:
\begin{lem}
\label{LC-1}
Let $\epsilon>0$. Let $f_\epsilon$ be a global weak solution to (L) with initial data satisfying the assumptions (\ref{assu1}-\ref{assu4}). Denote by $\rho_\epsilon(t,x):=\int f_\epsilon(t,x,v) dv$ and $J_\epsilon:= \int f_\epsilon v dv$. Then the following conservation laws hold in the distributional sense:
\begin{equation}
\partial_t \rho_\epsilon + \nabla_x. J_\epsilon=0. 
\end{equation}
\begin{equation}
\begin{split}
\partial_t J_\epsilon +  \nabla_x : \left(\int v \otimes v f_\epsilon dv \right) =&  -\frac{1}{2} \nabla_x (V_\epsilon+1)^2 \\
 +& \epsilon \operatorname{div}_x(\nabla_x V_\epsilon \otimes \nabla_x V_\epsilon) - \frac{\epsilon}{2}  \nabla_x \vert \nabla_x V_\epsilon \vert^2.
\end{split}
\end{equation}
\end{lem}

This is also the case for (S) and (S').
Indeed, the only difference with Lions-Perthame's equations is the nonlinear term $de^{V_\epsilon}$ in the Poisson equation, that we may consider as a source term by simply noticing that $de^{V_\epsilon} \in  L^1 \cap L^\infty(\mathbb{R}^3) $. More precisely we get, for any $\epsilon >0$, the (non-uniform) regularity estimates:

\begin{eqnarray}
\forall t \geq 0,  \int \vert v \vert^m f_{\epsilon}dvdx <\infty, \\
\rho_\epsilon \in L^\infty_{t, loc}(L^1_x\cap L^k_x).  
\end{eqnarray}

with $k=1+\frac{m}{3}$.

Thanks to these properties, we can prove that  the two first local conservation laws (\ref{LC1}) (conservation of charge)  and (\ref{LC2}) or (\ref{LC2'})  (conservation of momentum)  hold in the sense of distributions. 

\begin{lem}
\label{LC-2}
Let $\epsilon>0$. Let $f_\epsilon$ be a global weak solution to (S) or (S') with initial data satisfying the assumptions (\ref{assu1}-\ref{assu4}). We denote the two first moments  by $\rho_\epsilon(t,x):=\int f_\epsilon(t,x,v) dv$ and $J_\epsilon:= \int f_\epsilon v dv$ for any solution $f_\epsilon$ to (S) or (S').
The local conservation of charge reads:
\begin{equation}
\label{LC1}
\partial_t \rho_\epsilon + \nabla_x. J_\epsilon=0.
\end{equation}

The local conservation of current reads in the case of (S):
\begin{equation}
\label{LC2}
\partial_t J_\epsilon +  \nabla_x : \left(\int v \otimes v f_\epsilon dv \right) =   -d \nabla_x (e^{V_\epsilon}) + \epsilon \operatorname{div}_x(\nabla_x V_\epsilon \otimes \nabla_x V_\epsilon) - \frac{\epsilon}{2}  \nabla_x \vert \nabla_x V_\epsilon \vert^2,
\end{equation}
and in the case of (S'):
\begin{equation}
\label{LC2'}
\partial_t J_\epsilon +  \nabla_x : \left(\int v \otimes v f_\epsilon dv \right) =   -\frac{d}{\int d e^{V_\epsilon} dx} \nabla_x(e^{V_\epsilon}) + \epsilon \operatorname{div}_x(\nabla_x V_\epsilon \otimes \nabla_x V_\epsilon) - \frac{\epsilon}{2}  \nabla_x \vert \nabla_x V_\epsilon \vert^2.
\end{equation}
\end{lem}

\section{From (S) to Isothermal Euler}
\label{quasi-S}

The isothermal Euler equations (\ref{euler}) are hyperbolic symmetrizable. We can perform the change of unknown functions $(\rho,u) \mapsto (\log\frac{\rho}{d}, u)$ that leads to the system:
\begin{equation}
\label{neweuler}
\left\{
\begin{array}{ll}
\partial_t \log(\frac{\rho}{d}) + \nabla_x.  u + u.\nabla_x \log(\frac{\rho}{d}) - \nabla_x H.u=0, \\
\partial_t u +  u.\nabla_x u +\nabla_x \log(\frac{\rho}{d})=0.
\end{array}
\right.
\end{equation}
(we recall that by definition, $d= e^{-H}$).

Therefore, using classical results on hyperbolic symmetrizable systems (\cite{Maj}), we get the local  existence of smooth solutions:

\begin{prop}
\label{exi-iso} For any initial data $\rho_0>0, u_0$  such that $\rho_0 \in L^1(\mathbb{R}^3)$, $\log(\frac{\rho_0}{d}) \in H^s(\mathbb{R}^3)$ and $u_0 \in H^s(\mathbb{R}^3)$ for $s>\frac{3}{2}+1$, there is existence and uniqueness of a local smooth solution $\rho>0$ and $u$ to (\ref{neweuler}) such that :
\begin{equation}\label{reg1}\log\frac{\rho}{d}, u \in \mathcal{C}^0_t([0,T^*[, H^s(\mathbb{R}^3)) \cap \mathcal{C}^1_t([0,T^*[, H^{s-1}(\mathbb{R}^3))  \end{equation}
\begin{equation}\label{reg2}\rho \in \mathcal{C}^1([0,T^*[\times \mathbb{R}^3) \end{equation}
for some $T^*>0$.
\end{prop}

Since shocks may occur for large times, we will have to restrict to local results. We now prove the convergence of the charge and current density to the smooth solution to (\ref{neweuler}), as long as the latter exist.

We restrict to well-prepared, quasi monokinetic data.

\begin{thm}[The case of (S)] 
\label{theoS}Let $\rho_0> 0, u_0$ verifying the assumptions of Proposition \ref{exi-iso} and $\rho, u$ the corresponding strong solutions of system (\ref{euler}).
We assume that the sequence of initial data $(f_{0,\epsilon})$ satisfies the assumptions (\ref{assu1}-\ref{assu4}) and:
\begin{equation}
\label{condi2.1}
\int f_{0,\epsilon} \vert v -u_0\vert^2 dv dx \rightarrow 0,
\end{equation}
\begin{equation}
\label{condi2.2}
\Vert \sqrt{\epsilon} \nabla_x  V_{0,\epsilon}\Vert_{L^2}\rightarrow 0,
\end{equation}
\begin{equation}
\label{condi2.3}
\int \left(-de^{V_{0,\epsilon}} \log (de^{V_{0,\epsilon}}/\rho_{0}) -de^{V_{0,\epsilon}} + \rho_0 \right) dx \rightarrow 0,
\end{equation}
where $V_{0,\epsilon}$ is solution of the nonlinear Poisson equation:
$$
-\epsilon \Delta_x V_{0,\epsilon} = \int f_{0,\epsilon} dv - d e^{V_{0,\epsilon}}.
$$

Then $\rho_\epsilon$ weakly-* converges to $\rho$ and $J_\epsilon$ weakly-* converges to $\rho u$ in the weak sense of measures. Furthermore we have the following local in time strong convergences: $u_\epsilon=J_\epsilon/\rho_\epsilon$ strongly converges to $u$ in the following sense:
$$ \int {\vert u_\epsilon - u \vert^2}{\rho_\epsilon}dx \rightarrow 0$$
in $L^\infty_t$ and 
$$\sqrt{de^{V_\epsilon}} \rightarrow \sqrt{\rho}$$
in $L^\infty_t(L^2_x)$.
\end{thm}

\begin{rque}[On the class of admissible initial data satisfying assumptions (\ref{assu1}-\ref{assu4})and (\ref{condi2.1}-\ref{condi2.3})]
\label{coldions}
This class is not empty : indeed it includes Maxwellians of the form
\begin{equation}
f_{0,\epsilon} (x,v)=\frac{\rho_{0,\epsilon}(x)}{(2\pi T_{i,\epsilon})^{3/2}} e^{\frac{-\vert v -u_0(x)\vert^2}{2T_{i,\epsilon}}},
\end{equation}
where $\rho_{0,\epsilon}$ is computed by the Poisson equation after having previously chosen $V_{0,\epsilon}$ such that (\ref{condi2.2}) and (\ref{condi2.3}) hold ( for example, we can simply take $V_{0,\epsilon}=V_0$  with $V_0$ safisfying $de^{V_{0}}=\rho_0$ and (\ref{condi2.2}) and  (\ref{condi2.3}) trivially hold) and $T_{i,\epsilon} \rightarrow_{\epsilon \rightarrow 0} 0$ (cold ions approximation).

\end{rque}


\begin{proof} 

Let $(\rho,u)$ verifying the regularity of (\ref{reg1}-\ref{reg2}) (for the moment $\rho$ and $u$ do not a priori satisfy the isothermal Euler equations).

We recall that the energy for system (S) is the following functional:
 \begin{equation*}
\mathcal{F}_\epsilon(t)= \frac{1}{2}\int f_\epsilon\vert v \vert ^2 dv dx   + \frac{\epsilon}{2} \int \vert \nabla_x V_\epsilon \vert^2 dx  + \int d(V_\epsilon-1)e^{V_\epsilon} dx.
\end{equation*}

The first two terms correspond to the  energy for the quasineutral Vlasov-Poisson limit studied by Brenier in \cite{Br}. Therefore we accordingly modulate this quantity by considering:

$$\frac{1}{2}\int f_\epsilon\vert v - u\vert ^2 dv dx+ \frac{\epsilon}{2} \int \vert \nabla_x V_\epsilon \vert^2 dx.$$

Let us now look at $m_\epsilon:=de^{V_\epsilon}$. We can notice that:
$$\int d(V_\epsilon - 1)e^{V_\epsilon} dx= \int (m_\epsilon \log \left( m_\epsilon/ d \right)  - m_\epsilon)dx.$$

As mentioned in the introduction, we observe a strong analogy with the relative entropy in collisional kinetic equations (we refer to \cite{SR} for a reference on the topic). So by analogy, we modulate this quantity  and hence consider the following relative entropy:
\begin{equation}
\mathcal{H}_\epsilon(t)=\frac{1}{2}\int f_\epsilon\vert v - u\vert ^2 dv dx + \int (m_\epsilon \log \left( m_\epsilon/ \rho \right)  - m_\epsilon +\rho)dx + \frac{\epsilon}{2} \int \vert \nabla_x V_\epsilon \vert^2 dx.
\end{equation}

Later on, the well known inequality (which is a plain consequence of the inequality $x-1 \geq \log x $, for $x>0$) will be very useful:

\begin{equation}
\label{entropy}
\int (\sqrt{a}-\sqrt{b})^2 dx \leq \int (a \log(a/b) - a +b) dx.
\end{equation}



We want to show that  $\mathcal{H}_\epsilon(t)$ satisfies the inequality:
\begin{eqnarray}
{\mathcal{H}}_\epsilon(t) \leq   {\mathcal{H}}_\epsilon(0) +  G_\epsilon(t) +  C\int_0^t \Vert \nabla_x u \Vert_{L^\infty} \mathcal{H}_\epsilon(s) ds,
\end{eqnarray}
with  $G_\epsilon(t)\rightarrow 0 $ uniformly in time.

We show how we can deduce this kind of estimate (the computations can be rigorously justified using the local conservation laws of Lemma \ref{LC-2}). Since the energy is non-increasing, we have:
\begin{equation}
\frac{d{\mathcal{H}}_\epsilon(t) }{dt} \leq I_\epsilon(t),
\end{equation}	
with: 
\begin{equation}
\label{calcul2}
\begin{split}
I_\epsilon(t):= \int \partial_t f_\epsilon \left(\frac{1}{2} \vert u\vert^2 - v.u\right) dvdx +& \int f_\epsilon \partial_t \left(\frac{1}{2} \vert u\vert^2 - v.u\right) dv dx   \\   
+&  \int \partial_t\left(m_\epsilon \log(d/\rho)\right)dx + \int \partial_t \rho dx.
\end{split}
\end{equation}	

Let us first focus on the first two terms of $I_\epsilon(t)$. Thanks to the Vlasov equation satisfied by $f_\epsilon$ and after integrating by parts, we get:
\begin{equation}
\label{calcul3} 
\begin{split}
\int \partial_t f_\epsilon \left(\frac{1}{2} \vert u\vert^2 - v.u\right) dvdx +& \int f_\epsilon \partial_t \left(\frac{1}{2} \vert u\vert^2 - v.u\right) dv dx\\
=& \int f_\epsilon \left( \partial_t + v.\nabla_x + E_\epsilon. \nabla_v \right)\left(\frac{1}{2} \vert u\vert^2 - v.u\right) dvdx\\
=& \int f_\epsilon (u-v).(\partial_t + v.\nabla_x)u dvdx - \int f_\epsilon E_\epsilon. u dv dx\\
=& \int f_\epsilon (u-v).(\partial_t + u.\nabla_x)u dvdx -  \int f_\epsilon (u-v). \left( (u-v). \nabla_x u \right) dv dx\\
-&  \int \rho_\epsilon E_\epsilon. u dx .
\end{split}
\end{equation}

We can use the Poisson equation to compute the last term of (\ref{calcul3}).


\begin{eqnarray*}
-  \int \rho_\epsilon E_\epsilon. u dx &=&  \int \rho_\epsilon \nabla_x V_\epsilon dx \\
&=& \int de^{V_\epsilon}\nabla_x V_\epsilon u dx - \epsilon \int \Delta_x V_\epsilon \nabla_{x}V_\epsilon. udx \\
&=& \int  d \nabla_x e^{V_\epsilon}. u dx - \epsilon \int \nabla_x : (\nabla_x V_\epsilon \otimes \nabla_x V_\epsilon) u dx +\epsilon \int \frac{1}{2} \nabla_{x}\vert\nabla_x V_\epsilon\vert^2 u dx
\\
&=& -\int  de^{V_\epsilon} \operatorname{div}_x u dx -  \int e^{V_\epsilon} \nabla_x d .  u dx + \epsilon \int D(u)  : (\nabla_x V_\epsilon \otimes \nabla_x V_\epsilon)  dx\\
&-&\epsilon \int \frac{1}{2} \vert\nabla_x V_\epsilon\vert^2 \operatorname{div}_{x} u dx,  \\
\end{eqnarray*}
where $D(u) = \frac{1}{2} \left( \partial_{x_i} u_j + \partial_{x_j} u_i  \right)_{i,j}$ is the symmetric part of $\nabla_x u = (\partial_{x_i} u_j)_{i,j}$.

We now focus on the last two terms of $I_\epsilon(t)$:
\begin{equation}
\label{calcul4}
\begin{split}
\int \partial_t\left(m_\epsilon \log(d/\rho)\right)dx + \int \partial_t \rho dx
=&  \int (-de^{V_\epsilon}/\rho +1) \partial_t \rho + \int d \partial_t e^{V_\epsilon} \log(d/\rho) dx \\
=&  \int (-de^{V_\epsilon}/\rho +1) \partial_t \rho + \epsilon \int  \partial_t \Delta_{x} V_\epsilon \log(d/\rho) dx \\
-& \int \operatorname{div}_x J_\epsilon  \log(d/\rho) dx,\\
\end{split}
\end{equation}
using the Poisson equation and the local conservation of mass:
\[
\partial_t \rho_\epsilon = - \operatorname{div}_x J_\epsilon.
\]


We observe that:
\begin{equation}
\label{observation}
\begin{split}
\int f_\epsilon (u-v)& .\nabla_x \log\frac{\rho}{d} dvdx \\=&   \int \rho_\epsilon u.  \nabla_x\log \frac \rho d  dx - \int J_\epsilon. \nabla_x \log \frac \rho d dx \\
=&  \int de^{V_\epsilon}  u.  \nabla_x \log \frac \rho d  dx - \epsilon \int \Delta_{x}V_\epsilon u. \nabla_x \log \frac \rho d  dx  \\
+& \int J_\epsilon.  \nabla_x\log \frac \rho d dx.
\end{split}
\end{equation}

Consequently, according to (\ref{calcul3}) and (\ref{calcul4}) we get:

\begin{equation}
\begin{split}
I_\epsilon(t)=& \int -d e^{V_\epsilon}(\partial_t \rho +\nabla_x. u + u.\nabla_x \log \rho )dx+ \int \partial_t \rho dx\\
+& \int f_\epsilon(u-v) \left(\partial_t u + u.\nabla_x u + \nabla_x \log \frac \rho d\right)dvdx \\
+&\epsilon \int \Delta_{x}V_\epsilon u. \nabla_x \log \frac \rho d  dx + \epsilon \int  \partial_t \Delta_x V_\epsilon \log(d/\rho) dx  \\
-&  \int f_\epsilon (u-v) \left( (u-v). \nabla_x u \right) dv dx +\epsilon \int D(u)  : (\nabla_x V_\epsilon \otimes \nabla_x V_\epsilon)  dx \\
-&\epsilon \int \frac{1}{2} \vert\nabla_x V_\epsilon\vert^2 \operatorname{div}_{x} u dx.\\
\end{split}
\end{equation}

Let us now introduce the so-called acceleration operator $A$:

\begin{equation}
A(u,\rho)=\begin{pmatrix}\partial_t \log \frac \rho d +\nabla_x. u + u.\nabla_x \log \frac \rho d -\nabla_x H. u \\ \partial_t u + u.\nabla_x u +\nabla_x \log \frac \rho d
\end{pmatrix}.
\end{equation}

We observe that :
\begin{equation}
\begin{split}
\int \rho (\nabla_x.u  + u \nabla_x \log \rho) dx =& \int (\rho \nabla_x u + u.\nabla_x \rho)dx \\
=& \int \nabla_x.(\rho u)dx = 0,\\
\end{split}
\end{equation}
and thus $\int A(u,\rho).\begin{pmatrix}\rho \\ 0 \end{pmatrix} dx= \int \partial_t \rho dx$.

Gathering the pieces together we have proved:
\begin{equation}
\label{crucial}
\begin{split}
I_\epsilon(t)=& \int A(u,\rho).\begin{pmatrix}
-de^{V_\epsilon}+\rho  \\ \rho_\epsilon u-J_\epsilon\end{pmatrix}dx \\
+&\epsilon \int \Delta_{x}V_\epsilon u. \nabla_x \log \frac \rho d  dx + \epsilon \int  \partial_t \Delta_x V_\epsilon \log(d/\rho) dx  \\
-&  \int f_\epsilon (u-v) \left( (u-v). \nabla_x u \right) dv dx +\epsilon \int D(u)  : (\nabla_x V_\epsilon \otimes \nabla_x V_\epsilon)  dx \\
-& \epsilon \int \frac{1}{2} \vert\nabla_x V_\epsilon\vert^2 \operatorname{div}_{x} u dvdx.
\end{split}
\end{equation}

We are now ready to prove that $\mathcal{H}_\epsilon$ satisfies the expected stability inequality.

It is readily seen that there exists a constant independent of $\epsilon$ such that:
\begin{eqnarray*}
\left\vert \int f_\epsilon (u-v) \left( (u-v). \nabla_x u \right) dv dx\right\vert \leq C \int f_\epsilon \vert u-v \vert^2  \Vert \nabla_x u \Vert_{L^\infty} dv dx,\\
\left\vert \epsilon \int D(u)  : (\nabla_x V_\epsilon \otimes \nabla_x V_\epsilon)  dx \right\vert \leq C \int \epsilon \vert \nabla_x V_\epsilon \vert^2 \Vert \nabla_x u \Vert_{L^\infty} dx.
\end{eqnarray*}


We consider now:
\begin{eqnarray*}
G_\epsilon(t)&:=&\int_0^t  \left(\epsilon \int \Delta_{x}V_\epsilon u. \nabla_x \log \frac \rho d  dx + \epsilon \int  \partial_s \Delta_{x} V_\epsilon \log(d/\rho) dx \right)ds \\
&=& \int_0^t  \left(-\epsilon \int \nabla_{x}V_\epsilon. \nabla_x\left( u. \nabla_x \log \frac \rho d \right) dx - \epsilon \int \partial_s \nabla_{x} V_\epsilon .\nabla_x \log(d/\rho) dx \right)ds \\
&=&  \int_0^t   -\sqrt{\epsilon} \int \sqrt{\epsilon} \nabla_{x}V_\epsilon .\nabla_x \left( u.\nabla_x \log\left(\frac \rho d \right)\right) dx ds + \sqrt{\epsilon}\int_0^t \int \sqrt{\epsilon} \nabla_{x} V_\epsilon.  \partial_s \nabla_x \log(d/\rho) dx ds  \\
&-&  \sqrt{\epsilon}  \int  \sqrt{\epsilon}\nabla_{x}V_\epsilon(t,x) . \nabla_x \log(d/\rho(t,x)) dx + \sqrt{\epsilon}  \int \sqrt{\epsilon}\nabla_{x}V_\epsilon(0,x).   \nabla_x \log(d/\rho(0,x)) dx.
\end{eqnarray*}

Thanks to the conservation of the energy, $\sqrt{\epsilon} \nabla_x V_\epsilon$ is bounded uniformly with respect to $\epsilon$ in $L^\infty_t(L^2_x)$. Consequently, using Cauchy-Schwarz inequality, we get for any $0\leq t \leq T$:

\begin{eqnarray*}
G_\epsilon(t)&\leq& C \sqrt{\epsilon}\Vert \sqrt\epsilon \nabla_x V_\epsilon \Vert_{L^\infty_t(L^2_x)} \times\\
& &\Big( \Vert \nabla_x(u.\nabla_x\log{\rho/d}) \Vert_{L^\infty_t(L^2_x)}+ \Vert \log(\rho/d)   \Vert_{W^{1,\infty}_t(H^1_x)} \Big),
\end{eqnarray*}

and so we have $G_\epsilon(t) \rightarrow 0$ when $\epsilon \rightarrow 0$, locally uniformly in time.

Finally since
$$
{\mathcal{H}}_\epsilon(t) \leq {\mathcal{H}}_\epsilon(0) + \int_0^t I_\epsilon(s) ds,
$$
we have proved that:
\begin{eqnarray}
{\mathcal{H}}_\epsilon(t) \leq {\mathcal{H}}_\epsilon(0) + \int_0^t \int A(u,\rho).\begin{pmatrix}
-de^{V_\epsilon}+\rho  \\ \rho_\epsilon u - J_\epsilon \end{pmatrix}dx + G_\epsilon(t)  \nonumber \\ 
\label{stability} + C \left(\int_0^t \int f_\epsilon \vert u-v\vert^2 \Vert\nabla_x u \Vert_{L^\infty} dv dx ds +  \epsilon \int_0^t \int \frac{1}{2} \vert\nabla_x V_\epsilon\vert^2\Vert\nabla_x u \Vert_{L^\infty}  dx  ds\right).
\end{eqnarray}

Now we can choose $\rho$ and $u$ to be solutions of $A(\rho,u)=0$, with initial conditions $(\rho,u)_{\vert(t=0)}=(\rho_0,u_0)$. In other words $\rho$ and $u$ are solutions to the Isothermal Euler system (\ref{neweuler}).

Then we have:
\begin{eqnarray}
{\mathcal{H}}_\epsilon(t) \leq   {\mathcal{H}}_\epsilon(0) +  G_\epsilon(t) +  C\int_0^t \Vert \nabla_x u \Vert_{L^\infty} \mathcal{H}_\epsilon(s) ds
\end{eqnarray}
and thus, since $ {\mathcal{H}}_\epsilon(0)\rightarrow 0 $ and $G_\epsilon(t)\rightarrow 0 $, we deduce by Gronwall inequality that 

 \begin{equation}{\mathcal{H}}_\epsilon(t) \rightarrow 0,
 \end{equation}
 when $\epsilon \rightarrow 0$ (uniformly with respect to time).

By inequality (\ref{entropy}), this means in particular that 
$$\sqrt{de^{V_\epsilon}} \rightarrow \sqrt{\rho}.$$ 
strongly in $L^\infty_tL^2_x$.

Because of the uniform estimates in $L^\infty_t(L^1_x)$, $\rho_\epsilon$ (resp. $J_\epsilon$) weakly converges (up to a subsequence) in the weak sense of measures to some $\tilde{\rho}$ (resp. $J$).
In the other hand, in the sense of distributions, thanks to the quasineutral Poisson equation:
$$de^{V_\epsilon}-\rho_\epsilon \rightharpoonup 0$$
and thus,
$${de^{V_\epsilon}} \rightharpoonup\tilde{ \rho}.$$ 

Therefore, by uniqueness of the limit we deduce that $\tilde{\rho}= \rho$.

The last step of the proof relies on a by now classical convexity argument. We first get the following Cauchy-Schwarz inequality:

\begin{equation}
\frac{\vert J_\epsilon - \rho_\epsilon u\vert^2}{\rho_\epsilon}= \frac{\left(\int f_\epsilon(v-u)dv\right)^2}{\int f_\epsilon dv}\leq  \int f_\epsilon\vert v-u\vert^2 dv.
\end{equation}

The functional $(\rho,J) \rightarrow \int \frac{\vert J - \rho u\vert^2}{\rho} dx$ is convex and lower semi-continuous with respect to the weak convergence of measures (see \cite{Br}). Consequently the weak convergence in the sense of measures $\rho_\epsilon \rightharpoonup \rho$ and $J_\epsilon \rightharpoonup J$ leads to:

\begin{equation}
\int  \frac{\vert J - \rho u \vert^2}{\rho} dx \leq \liminf_{\epsilon\rightarrow 0} \int \frac{\vert J_\epsilon - \rho_\epsilon u\vert^2}{\rho_\epsilon}dx.
\end{equation}
So $J=\rho u$.

To conclude, the uniqueness of the limit allows us to say that the weak convergences actually hold without any extraction.

\end{proof}

\begin{rque}[Rate of convergence]Assume that $\mathcal{H}_{\epsilon}(0) \leq C \sqrt{\epsilon}$. Then the previous estimates show that locally uniformly in time:
\begin{equation}
\mathcal{H}_{\epsilon}(t) \leq C \sqrt{\epsilon}.
\end{equation}
\end{rque}

\section{Generalization to other quasineutral limits}
\label{quasi-others}

\subsection{From (S') to Isothermal Euler}

Similarly, we can prove an analogous theorem for system (S'):

\begin{thm}[The case of (S')] 
\label{theoS'}Let $\rho_0> 0, u_0$ verifying the assumptions of Proposition \ref{exi-iso}  and $\rho, u$ the corresponding strong solutions of system (\ref{euler}).
We assume that the sequence of initial data $(f_{\epsilon,0})$ satisfies the assumptions (\ref{assu1}-\ref{assu4}) and:
\begin{equation}
\label{condi2.1'}
\int f_{0,\epsilon} \vert v -u_0\vert^2 dv dx \rightarrow 0,
\end{equation}

\begin{equation}
\label{condi2.2'}
\sqrt{\epsilon} \nabla_x  V_{0,\epsilon}\rightarrow 0,
\end{equation}
strongly in $L^2$ 
and
\begin{equation}
\label{condi2.3'}
\int \left(-\frac{de^{V_{0,\epsilon}}}{\int de^{V_{0,\epsilon}}dx} \log \left(\frac{de^{V_{0,\epsilon}}}{\int de^{V_{0,\epsilon}}dx}/\rho_{0}\right) - \frac{de^{V_{0,\epsilon}}}{\int de^{V_{0,\epsilon}}dx} +\rho_0 \right) dx \rightarrow 0,
\end{equation}
where $V_{0,\epsilon}$ is solution of the nonlinear Poisson equation:
$$
-\epsilon \Delta_x V_{0,\epsilon} = \int f_{0,\epsilon} dv -  \frac{de^{V_{0,\epsilon}}}{ \int d e^{V_{0,\epsilon}} dx}.
$$

Then $\rho_\epsilon$ weakly-* converges to $\rho$ and $J_\epsilon$ weakly-* converges to $\rho u$ in the weak sense of measures. Furthermore, we have the following local strong convergences: $u_\epsilon=J_\epsilon/\rho_\epsilon$ strongly converges to $\rho u$ in the following sense:
$$ \int{\vert u_\epsilon -  u \vert^2}{\rho_\epsilon}dx \rightarrow 0$$

in $L^\infty_t$ and 
$$\sqrt{\frac{de^{V_\epsilon}}{\int d e^{V_\epsilon} dx}} \rightarrow \sqrt{\rho}$$
in $L^\infty_t(L^2_x) $.

\end{thm}

\begin{proof}[Sketch of proof] According to Theorem \ref{exiS}, the functional $\mathcal{G}_\epsilon(t)$ is non-increasing:
\begin{eqnarray*}
\mathcal{G}_\epsilon(t)= \frac{1}{2}\int f_\epsilon\vert v \vert ^2 dv dx + \int d(x)\left(V_\epsilon - \log \left(\int d e^{V_\epsilon} dx\right)\right)\frac{e^{V_\epsilon}}{\int d e^{V_\epsilon} dx} dx \nonumber \\
   + \frac{\epsilon}{2} \int \vert \nabla_x V_\epsilon \vert^2 dx .
\end{eqnarray*}

As in the previous proof, we consider $m_\epsilon = \frac{de^{V_\epsilon}}{\int d e^{V_\epsilon} dx}$, and we notice that:

\begin{equation}
\int d\left(V_\epsilon - \log \left(\int d e^{V_\epsilon} dx\right)\right)\frac{e^{V_\epsilon}}{\int d e^{V_\epsilon} dx} dx = \int (m_\epsilon \log(m_\epsilon/d))dx.
\end{equation}

Since $\int m_\epsilon dx=1$ (and thus, $\partial_t \int m_\epsilon dx=0$) we can actually add this quantity to the energy so that:
\begin{equation}
\int d\left(V_\epsilon - \log \left(\int d e^{V_\epsilon} dx\right)-1\right)\frac{e^{V_\epsilon}}{\int d e^{V_\epsilon} dx} dx = \int (m_\epsilon \log(m_\epsilon/d) -m_\epsilon) dx,
\end{equation}
and we are in the same case as before. Therefore we can consider the same modulated energy $\mathcal{H}_\epsilon(t)$. We skip the computations, which are very similar.

\end{proof}

\subsection{From (L) to Shallow-Water}

We now treat the case of system (L), following the same methodology as before.

\subsubsection{Formal derivation of the Shallow-Water equations}


For monokinetic data, i.e. $f_\epsilon(t,x,v) = \rho_\epsilon (t,x) \delta_{v=u_\epsilon(t,x)}$, the conservation laws state:

\begin{equation}
\partial_t \rho_\epsilon + \nabla_x. (\rho_\epsilon u_\epsilon)=0,
\end{equation}

\begin{equation}
\partial_t (\rho_\epsilon u_\epsilon) +  \nabla_x : (\rho_\epsilon u_\epsilon \otimes u_\epsilon) = -\rho_\epsilon \nabla_x V_\epsilon.
\end{equation}

We recall that the Poisson equation is:
\begin{equation}
V_\epsilon -\epsilon \Delta_x V_\epsilon = \int f_\epsilon dv -1.
\end{equation}

Since $\rho_\epsilon$ and $J_\epsilon$ are uniformly bounded in $L^\infty_t(L^1_x)$, the following convergences (up to a subsequence) hold in the sense of distributions: $\rho_\epsilon \rightharpoonup \rho$ and $J_\epsilon \rightharpoonup J$.

If we formally pass to the limit $\epsilon \rightarrow 0$ we get:
\begin{equation}
V= \rho-1.
\end{equation}

Consequently the limit system is the following:
\begin{equation}
\left\{
\begin{array}{ll}
\partial_t \rho + \nabla_x. (\rho u)=0, \\
\partial_t (\rho u) +  \nabla_x : (\rho u \otimes u) = - \rho \nabla_x \rho.
\end{array}
\right.
\end{equation}
or equivalently for smooth data:

\begin{equation}
\label{sw}
\left\{
\begin{array}{ll}
\partial_t \rho + \nabla_x. (\rho u)=0, \\
\partial_t u +  u.\nabla_x u = - \nabla_x \rho.
\end{array}
\right.
\end{equation}
As it has been said before, this system can be interpreted in $1D$ or $2D$ as the Shallow Water equations.

\begin{rque}[On the kinetic version of the shallow-water limit]
In a formal sense, one can also easily perform the kinetic limit $\epsilon \rightarrow 0$ and get the equation:
\begin{equation}
\left\{
\begin{array}{ll}
\partial_t f + v.\nabla_x f - \nabla_x  \rho. \nabla_v f=0, \\
\rho=\int f dv.
\end{array}
\right.
\end{equation}

To our knowledge, this equation is very badly mathematically understood.
The only existence result we are able to prove is the local existence of analytic solutions. Actually the proof given by Mouhot and Villani \cite{MV} (section $9$, local in time interaction) in the Vlasov-Poisson case identically holds. Indeed, we notice that in that proof, they do not need the smoothing effect on the force field provided by the Poisson equation. Although it is not explicitly said, the case of the singular force field $F= \nabla_x \rho$ is automatically included in their analysis. Of course this is not the case for the other results of their paper.
\end{rque}

\begin{rque}
Let us also mention that the quasineutral together with the gyrokinetic limit of a similar system was performed by the author in \cite{DHK1}. With quite general initial data, we get a limit equation of kinetic nature. In other words we do not need to restrict to particular initial data (or to strong regularity); this rather remarkable fact is due to the anisotropy of the system with the so-called finite Larmor radius scaling (\cite{FS2}). The Poisson equation then only degenerates in the magnetic field direction but this is overcome thanks to an averaging lemma.
\end{rque}

\subsubsection{Rigorous derivation for (partially) well-prepared data}

As  the Shallow Water equations (\ref{sw}) are hyperbolic symmetrizable we get the local existence of smooth solutions \cite{Maj}.

\begin{prop}For any initial data $\rho_0, u_0$ in $H^s(\mathbb{T}^n)$ for $s>\frac{n}{2}+1$, there is existence and uniqueness of a local smooth solution to (\ref{sw}):
$$\rho, u \in C^0_t([0,T^*[, H^s(\mathbb{T}^n)) \cap C^1_t([0,T^*[, H^{s-1}(\mathbb{T}^n))$$
for some $T^*>0$.
\end{prop}

As before, we restrict to finite time intervals.




\begin{thm}[The case of (L)]
\label{theoL} Let $\rho_0\geq 0, u_0 \in H^s(\mathbb{T}^n)$ ($s>n/2 +1$ large enough) and $\rho, u$ the corresponding strong solutions to System \ref{sw}.
We assume that the sequence of initial data $(f_{0,\epsilon})$ satisfies the hypotheses (\ref{assu1}-\ref{assu4}) and:
\begin{equation}
\label{condi1.1}
\int f_{0,\epsilon} \vert v -u_0\vert^2 dv dx \rightarrow 0,
\end{equation}
\begin{equation}
\label{condi1.2}
\int \vert  \sqrt{\epsilon}\nabla_x V_{0,\epsilon} \vert ^2  dx \rightarrow 0,
\end{equation}
and
\begin{equation}
\label{condi1.3}
(Id-\epsilon\Delta_x)^{-1}(\rho_{0,\epsilon}-1) \rightarrow (\rho_{0}-1) ,
\end{equation}
strongly in $L^2$.

Then $\rho_\epsilon$ weakly-* converges to $\rho$ and $J_\epsilon$ weakly-* converges to $\rho u$ in the weak sense of measures. Furthermore, we have the local strong convergences:

$$V_\epsilon\rightarrow \rho-1$$
in $L^\infty_t L^2_x$, $u_\epsilon$ strongly converges to $u$ in the following sense:
$$ \int {\vert u_\epsilon -  u \vert^2}{\rho_\epsilon}dx \rightarrow 0.$$

Moreover,
$$
\int \vert  \sqrt{\epsilon}\nabla_x V_{\epsilon}\vert ^2  dx \rightarrow 0.
$$

\end{thm}

\begin{rque}
\begin{itemize}
 \item Assumptions (\ref{condi1.2}) and (\ref{condi1.3}) are satisfied for some smooth $\rho_0$ as soon as $\rho_{0,\epsilon}$ strongly converges to $\rho_0$ in $L^2$ (For instance, when $\rho_{0,\epsilon}$ is uniformly bounded in some $H^\alpha$ with $\alpha >0$).
 Indeed, in this case, we notice that:
\begin{equation}
\frac{1}{\sqrt{\epsilon}}\Big((Id-\epsilon\Delta_x)^{-1}\rho_{0,\epsilon} - \rho_{0,\epsilon} \Big)
\end{equation}
lies in a compact of $H^{-1}$ endowed with its strong topology.

This implies, thanks to the Poisson equation that $\sqrt{\epsilon} \Delta_x V_{0,\epsilon}$ lies in a compact for the $H^{-1}$ norm. This means that $\sqrt{\epsilon} \nabla_x V_{0,\epsilon}$ strongly converges to some $\nabla_x \Psi_0$ in the $L^2$ norm (up to a sequence). But $V_{\epsilon,0}$ also strongly converges in $L^2$ to $\rho_{0}-1$, so by uniqueness of the limit in the sense of distributions, $\nabla_x \Psi_0=0$.

\item For this reason,  an "ill-prepared" case would correspond to some oscillating in space initial data: 
$$V_{0,\epsilon}= \rho_0 -1  + r_\epsilon,$$ 
where $r_\epsilon$ only weakly converges to $0$ and $\sqrt{\epsilon}\nabla_x r_\epsilon$ is bounded in $L^2$ (one can think of $r_\epsilon(x) = e^{ix/\sqrt{\epsilon}}$). Then we would have to filter these oscillations in space to prove strong convergences.

This indicates that the massless electrons have a stabilizing effect on the system, insofar as no time oscillations occur unless some space oscillations are imposed at the initial time.

\end{itemize}

\end{rque}

\begin{rque}
We need some additional regularity on $\rho$ and $u$ in order to handle some non-linear quantities : so we take $s$ large enough (the lower bound $n/2 +1$ is not sufficient). But we will not dwell on the optimal constant.
\end{rque}

\begin{proof}[Sketch of proof] 
The functional $\mathcal{E}_\epsilon(t)$ recalled below is the energy of system (L):
$$
\mathcal{E}_\epsilon(t)= \frac{1}{2}\int f_\epsilon\vert v \vert ^2 dv dx + \frac{1}{2}\int V_\epsilon^2 dx  + \frac{\epsilon}{2} \int \vert \nabla_x V_\epsilon \vert^2 dx .
$$

The principle of the proof is to consider the following modulation of $\mathcal{E}_\epsilon(t)$:
\begin{equation}
\tilde{\mathcal{H}}_\epsilon(t)=\frac{1}{2}\int f_\epsilon\vert v - u\vert ^2 dv dx + \frac{1}{2}\int \vert V_\epsilon - ({\rho}-1) \vert^2 dx  + \frac{1}{2} \int \vert \sqrt{\epsilon} \nabla_x V_\epsilon \vert^2 dx .
\end{equation}

Then we can show with similar considerations as previous proofs that $\mathcal{H}_\epsilon(t)$ satisfies an inequality of the form:
 \begin{eqnarray*}
\tilde {\mathcal{H}}_\epsilon(t) \leq   \tilde{\mathcal{H}}_\epsilon(0)+ \tilde G_\epsilon(t) +  \int_0^t C\Vert \partial_x u \Vert_{L^\infty} \tilde{\mathcal{H}}_\epsilon(s)  ds
\end{eqnarray*}
with $\tilde G_\epsilon(t) \rightarrow0$ when $\epsilon \rightarrow 0$, uniformly in time.

\end{proof}

\subsection{Quasineutral limit for the isothermal Euler-Poisson system of Cordier and Grenier}
\label{IEP}

As it was mentioned in the introduction,  in \cite{CG}, Cordier and Grenier study  an isothermal Euler-Poisson version of (S2) and prove the quasineutral limit to the same kind of Euler equation (\ref{euler}). So our result can be seen somehow as a generalization of theirs, since our startpoint is the kinetic equation. Actually the relative entropy method can also apply to their system.

In \cite{CG}, Cordier and Grenier consider the isothermal Euler-Poisson system (in $1D$):
 \begin{equation}
 \label{IEPequ}
  \left\{
      \begin{array}{ll}
\partial_t \rho_\epsilon + \partial_x (\rho_\epsilon u_\epsilon)=0 \\
\partial_t u_\epsilon + u_\epsilon\partial_x u_\epsilon + \frac{T}{\rho_\epsilon} \partial_x \rho_\epsilon = - \partial_x V_\epsilon \\
-\epsilon \partial_{xx}^2 V_\epsilon = \rho_\epsilon - e^{V_\epsilon} ,
\end{array}
  \right.
\end{equation}
where $T$ is the (scaled) temperature of ions, of order $1$.

The authors perform the quasineutral limit $\epsilon \rightarrow 0$ to the so-called quasineutral Euler system:
 \begin{equation}
  \left\{
      \begin{array}{ll}
\partial_t \rho + \partial_x (\rho u)=0 \\
\partial_t u + u\partial_x u + \frac{T+1}{\rho} \partial_x \rho = 0. \\
\end{array}
  \right.
\end{equation}
Their proof relies on rather tricky energy estimates obtained by pseudodifferential calculus (using the framework introduced by Grenier \cite{Gre3}). 

Relative entropy methods provide an alternative and more direct proof. Indeed, the following functional is an energy for system (\ref{IEPequ}):
\begin{equation}
\mathfrak{E}_\epsilon(t) = \frac{1}{2} \int \rho_\epsilon u_\epsilon^2 dx + T \int \rho_\epsilon(\log \rho_\epsilon  - 1) dx  + \int (V_\epsilon-1)e^{V_\epsilon} dx  + \frac{\epsilon}{2} \int \vert \partial_x V_\epsilon \vert^2 dx.
\end{equation}

We can consequently consider the modulated energy:
\begin{equation}
\begin{split}
\mathfrak{H}_\epsilon(t) =& \frac{1}{2} \int \rho_\epsilon\vert  u_\epsilon - u \vert ^2 dx + T \int \rho_\epsilon(\log \left(\rho_\epsilon/\rho\right)  - 1 + \rho/\rho_\epsilon ) dx  \\+&\int (m_\epsilon \log \left( m_\epsilon/ \rho \right)  - m_\epsilon +\rho)dx  
+ \frac{\epsilon}{2} \int \vert \partial_x V_\epsilon \vert^2 dx,
\end{split}
\end{equation}
with $m_\epsilon = e^{V_\epsilon}$. We can show, as in the previous proofs, that $H_\epsilon$ satisfies some stability inequality.  In addition to the results obtained for system (S), we get the following strong convergence in $L^\infty_t(L^2_x)$:
\begin{equation}
\sqrt{\rho_\epsilon} \rightarrow \sqrt{\rho}.
\end{equation}

\section{Combined quasineutral and large magnetic field limit}
\label{quasi-magn}

We now study the limit $\epsilon \rightarrow 0$ of the following system, which is nothing but system (S') with a strong magnetic field.

\begin{equation}
\label{S'rot}
\left\{
    \begin{array}{ll}
  \partial_t f_\epsilon + v.\nabla_x f_\epsilon + \left(E_\epsilon+ \frac{v\wedge e_\parallel}{\epsilon}\right).\nabla_v f_\epsilon =0  \\
  E_\epsilon = -\nabla_x V_\epsilon\\
-\epsilon^{2\alpha} \Delta_x V_\epsilon = \int f_\epsilon dv -\frac{d e^{V_\epsilon}}{\int d e^{V_\epsilon} dx}\\
    f_{\epsilon,\vert t=0} =f_{0,\epsilon}, \int f_{0,\epsilon} dvdx =1.\\
    \end{array}
  \right.
\end{equation}

We shall not dwell on the existence of global weak solutions, since it is very similar to the theory for system (S') that was studied in section \ref{GW}.

We start with a formal analysis in order to show how we can get the expected limit system.

\subsection{Formal analysis}

For monokinetic data, i.e. $f_\epsilon(t,x,v) =\rho_\epsilon(t,x) \delta(v= u_\epsilon(t,x))$, the two first conservation laws read:
$$
\partial \rho_\epsilon + \nabla_x.(\rho_\epsilon u_\epsilon)=0,
$$
$$
\partial_t u_\epsilon + u_\epsilon. \nabla_x u_\epsilon = E_\epsilon  + \frac{u_\epsilon^\perp}{\epsilon}.
$$
The Poisson equation reads:

$$
-\epsilon^{2\alpha} \Delta_x V_\epsilon = \rho_\epsilon -\frac{d e^{V_\epsilon}}{\int d e^{V_\epsilon}dx}.
$$

In the limit $\epsilon\rightarrow 0$, assuming that in some sense $\rho_\epsilon \rightharpoonup \rho, V_\epsilon \rightharpoonup V$ 
(as well as $\nabla_x V_\epsilon \rightharpoonup \nabla_x V $), we get:
\begin{equation}
 \rho= \frac{d e^{V}}{\int d e^{V}dx},
\end{equation}
and this implies that:
$$
\nabla_x V = \frac{\nabla_x \rho}{\rho} -  \frac{\nabla_x d}{d}.
$$

If we multiply the second conservation law by $\epsilon$ we get:
$$
u_\epsilon^\perp = \epsilon \left(\partial_t u_\epsilon +u_\epsilon. \nabla_x u_\epsilon -E_\epsilon \right).
$$
This implies that $u_\epsilon^\perp \rightharpoonup 0$. This convergence can not occur in a strong sense because of the oscillations in time of frequency $\mathcal{O}(1/\epsilon)$, created by the magnetic field, but we can precisely describe the oscillations and consequently the strong convergence.

We denote by $\mathcal{R}(s)$ the rotation of axis $e_\parallel$ and angle $s$. Explicityly, we have:

$$
\mathcal{R}(s)=\begin{pmatrix} \cos t & -\sin t  & 0 \\ \sin t & \cos t & 0 \\ 0 & 0 & 1\end{pmatrix}.
$$

Following standard methods for singular perturbation problems (\cite{Gre5}, \cite{Sch}), we introduce the filtered momentum $w_\epsilon$ defined by:
\begin{equation}
w_\epsilon = \mathcal{R}(t/\epsilon) u_\epsilon.
\end{equation}

 Then it is readily seen that $w_\epsilon$ satisfies the equation:
 
\begin{equation}
\partial_t w_\epsilon + \mathcal{R}(-t/\epsilon)w_\epsilon. \nabla_x w_\epsilon =  \mathcal{R}(t/\epsilon) E_\epsilon.
\end{equation}

We assume then that $w_\epsilon \rightarrow w$ strongly. We take the limit
$\epsilon\rightarrow 0$ by time averaging:

\begin{equation}
\mathcal{R}(-t/\epsilon)w_\epsilon. \nabla_x w_\epsilon \rightarrow \frac{1}{2\pi} \int_0^{2\pi} \mathcal{R}(-\tau)w_\epsilon. \nabla_x w_\epsilon d\tau=w_\parallel \partial_{x_\parallel} w
\end{equation}
and argue similarly for the other terms.
We get in the end the following isothermal Euler system (with no dynamics in the $x_\perp$ variable):
\begin{equation}
\label{eulerrotating}
\left\{
    \begin{array}{ll}
\partial_t \rho + \partial_{x_\parallel} (\rho w_\parallel) =0, \\
\partial_t w + w_\parallel \partial_{x_\parallel} w = -\frac{\nabla_{x_\parallel} \rho}{\rho} - \nabla_{x_\parallel} H.
    \end{array}
  \right.
\end{equation}

Of course, this system is very similar to the "usual" isothermal Euler system (\ref{euler}), so we get the same existence result.

\begin{prop}
\label{exi-iso2} For any initial data $\rho_0>0, w_0$  such that $\rho_0 \in L^1(\mathbb{R}^3)$, $\log(\frac{\rho_0}{d}) \in H^s(\mathbb{R}^3)$ and $w_0 \in H^s(\mathbb{R}^3)$ for $s>\frac{3}{2}+1$, there is existence and uniqueness of a local smooth solution $\rho>0$ and $u$ to (\ref{neweuler}) such that :
\begin{equation}\log\frac{\rho}{d}, w \in \mathcal{C}^0_t([0,T^*[, H^s(\mathbb{R}^3)) \cap \mathcal{C}^1_t([0,T^*[, H^{s-1}(\mathbb{R}^3))  \end{equation}
\begin{equation}\rho \in \mathcal{C}^1([0,T^*[\times \mathbb{R}^3) \end{equation}
for some $T^*>0$.
\end{prop}

From what we have seen, we can thus expect the strong convergence
\begin{equation}
\mathcal{R}(t/\epsilon)u_\epsilon \rightarrow w,
 \end{equation}
that is to say:
\begin{equation}
u_\epsilon - \mathcal{R}(-t/\epsilon)w \rightarrow 0.
\end{equation}

\subsection{Convergence proof}

We first give the stability inequality we obtain for system (\ref{S'rot}).

\begin{prop}
Let $(f_{0,\epsilon})$ be a sequence of initial data satisfying assumptions  (\ref{assu1}-\ref{assu4}) and $(f_\epsilon)$ the corresponding global weak solutions to (\ref{S'rot}).

Let $s>3/2 +1$. For any sequence $\log \frac{\bar{\rho}_\epsilon}{d}$, $\bar{u}_\epsilon$ in $\mathcal{C}^0_t([0,T^*[, H^s(\mathbb{R}^3)) \cap \mathcal{C}^1_t([0,T^*[, H^{s-1}(\mathbb{R}^3))$  we define the modulated energy:
\begin{equation}
\mathcal{H}_\epsilon(t)=\frac{1}{2}\int f_\epsilon\vert v - \bar{u}_\epsilon\vert ^2 dv dx + \int (m_\epsilon \log \left( m_\epsilon/ \bar{\rho}_\epsilon \right)  - m_\epsilon +\bar{\rho}_\epsilon)dx + \frac{\epsilon^{2\alpha}}{2} \int \vert \nabla_x V_\epsilon \vert^2 dx,
\end{equation}
with $m_\epsilon = \frac{de^{V_\epsilon}}{\int d e^{V_\epsilon} dx} $.
Then the following inequality holds:
\begin{eqnarray}
\label{stability-B}
{\mathcal{H}}_\epsilon(t) \leq   {\mathcal{H}}_\epsilon(0) +  G_\epsilon(t) +  C\int_0^t \Vert \nabla_x \bar{u}_\epsilon \Vert_{L^\infty} \mathcal{H}_\epsilon(s) ds \nonumber \\
+ \int_0^t \int A_\epsilon(\bar \rho_\epsilon, \bar u_\epsilon). \begin{pmatrix} -m_\epsilon +\bar{\rho}_\epsilon \\ J_\epsilon -\rho_\epsilon\bar{u}_\epsilon\end{pmatrix} dx ds,
\end{eqnarray}
with  $A_\epsilon(t,x)$ the so-called acceleration operator defined by:
\begin{equation}
\label{acc-op}
A_\epsilon({\bar \rho_\epsilon}, \bar u_\epsilon)= \begin{pmatrix}\partial_t \log\left(\frac{\bar\rho_\epsilon}{d}\right) + \nabla_x.\bar u_\epsilon + \bar u_\epsilon. \nabla_x  \log\left(\frac{\bar\rho_\epsilon}{d}\right) - \nabla_x H . \bar u_\epsilon\\  \partial_t \bar u_\epsilon + \bar u_\epsilon. \nabla_x \bar u_\epsilon  + \nabla_x \log\left(\frac{\bar\rho_\epsilon}{d}\right)- \frac{\bar u_\epsilon^\perp}{\epsilon} \end{pmatrix},
\end{equation}
and $G_\epsilon(t)$ satisfying:
\begin{eqnarray}
G_\epsilon(t)&\leq& C{\epsilon}^\alpha\Vert\epsilon^\alpha \nabla_x V_\epsilon \Vert_{L^\infty_t(L^2_x)} \times \nonumber\\
& &\Big( \Vert \nabla_x(\bar u_\epsilon.\nabla_x\log{\bar\rho_\epsilon/d}) \Vert_{L^\infty_t(L^2_x)}+ \Vert \log(\bar\rho_\epsilon/d)   \Vert_{W^{1,\infty}_t(H^1_x)} \Big) \label{esti-G}.
\end{eqnarray}
\end{prop}

\begin{proof}
The proof is similar to the one given to obtain (\ref{stability}) in the proof of Theorem \ref{theoS} and therefore we omit it.
\end{proof}

As we wish to show the strong convergence (after filtering) of $(\rho_\epsilon:=\int f_\epsilon dv, u_\epsilon:= \frac{1}{\rho_\epsilon} \int v f_\epsilon dv)$ to solutions to system (\ref{eulerrotating}),
a natural idea would consist in taking $(\bar\rho_\epsilon:= \rho, \bar u_\epsilon= \mathcal{R}(-t/\epsilon)w)$ where $\rho$ and $w$ are the smooth solution to system (\ref{eulerrotating}) with initial data $\rho_0$ and $w_0$.
Unfortunately, we can not prove directly 
$$A_\epsilon(\rho, \mathcal{R}(-t/\epsilon)w) \rightarrow 0$$
in a strong sense. Thus, as in \cite{Sch} or \cite{GSR2}, we add a small correction denoted by $\epsilon z_\epsilon$ in order to build a higher order approximation of the equation and make the acceleration operator vanish. 
The shape of $z_\epsilon$ is precisely chosen in order to ``kill'' the non-vanishing terms which only weakly converge to $0$.

In the following, we will consider the derivative of  $\mathcal{R}$:
 
 \begin{equation}
 \label{defS}					
 \mathcal{S}(t):= \frac{d \mathcal{R}(t)}{dt},
\end{equation}
which satisfies: 
$$
 \frac{d \mathcal{S}(t)}{dt}= \mathcal{R}(-t)_\perp:=\begin{pmatrix} \cos t & \sin t  & 0 \\ -\sin t & \cos t & 0 \\ 0 & 0 & 0\end{pmatrix}.
$$


\begin{thm}
\label{theoS'rot} Let $\rho_0, w_0$ initial data verifying the hypotheses of Proposition \ref{exi-iso2} with $s> 5/2$. We assume that the sequence of initial data $(f_{0,\epsilon})$ satisfies the assumptions (\ref{assu1}-\ref{assu4}) and:
\begin{equation}
\mathcal{H}_\epsilon(0) \rightarrow 0.
\end{equation}
Let $(\log\frac{\rho}{d}, w)$ the unique strong solution to (\ref{eulerrotating}) with $\left(\log\frac{\rho_0}{d}, w_0\right)$ as initial conditions. We define $\bar\rho_\epsilon$ and $\bar u_\epsilon$ by the relation:
\begin{equation}
\label{defbar}
\begin{pmatrix} \log\frac{\bar \rho_\epsilon}{d} \\ \bar u_\epsilon \end{pmatrix}=  \begin{pmatrix}\log\frac{\rho}{d}\\ \mathcal{R}(-t/\epsilon) w \end{pmatrix} +  \epsilon y_\epsilon,
\end{equation}
with $y_\epsilon = \begin{pmatrix} z^{\rho}_\epsilon \\ \mathcal{R}(-t/\epsilon) z^{w}_\epsilon \end{pmatrix}$ and  $z_\epsilon^{\rho}$ (resp. $z_\epsilon^{\rho}$) defined by its Fourier transform $\mathcal{F}z_\epsilon^{\rho}$ (resp. $\mathcal{F}z_\epsilon^{w}$):


\begin{equation}
\begin{split}
\mathcal{F}z^\rho_\epsilon(\xi) =&-\mathbbm{1}_{\vert \xi \vert \leq\frac{1}{\epsilon}}\mathcal{F}\left(\nabla_{x_\perp}. (\mathcal{S}(t/\epsilon)w)  - \nabla_x H . \mathcal{S}(t/\epsilon) w_\perp \right) \\
-& \int_{\mathbb{R}^3} \mathbbm{1}_{\vert \xi- \eta \vert + \vert \eta \vert \leq\frac{1}{\epsilon}}\mathcal{F}(\mathcal{S}(t/\epsilon) w_\perp)(\eta). \mathcal{F} \left(\nabla_{x_\perp}  \log\left(\frac{\rho}{d}\right)\right)(\xi-\eta)d\eta,
\end{split}
\end{equation}

\begin{equation}
\begin{split}
\mathcal{F}z^{w}_\epsilon(\xi) =&  -\mathbbm{1}_{\vert \xi \vert \leq\frac{1}{\epsilon}}\mathcal{F}\left( \mathcal{S}(t/\epsilon)\nabla_{x_\perp}\log\left(\frac{\rho}{d}\right)\right) \\
-& \int_{\mathbb{R}^3} \mathbbm{1}_{\vert \xi- \eta \vert + \vert \eta \vert \leq\frac{1}{\epsilon}} \mathcal{F}(\mathcal{S}(t/\epsilon) w_\perp)(\eta).\mathcal{F}(\nabla_{x_\perp} w)(\xi-\eta) d\eta  , 
\end{split}
\end{equation}
where the operator $\mathcal{S}(t)$ is defined in (\ref{defS}).

Then, there exists $C>0$ depending only on $w$ and $\log \frac \rho d$, 
$$\Vert z_\epsilon \Vert_{L^\infty_t([0,T], H^{s-1})} \leq  C, $$
and locally  uniformly in time we have:
\begin{equation}
\mathcal{H}_\epsilon(t) \rightarrow 0. 
\end{equation}

In particular, this means that $\rho_\epsilon$ weakly-* converges to $\rho$ and $J_\epsilon$ weakly-* converges to $\rho w_\parallel$. Furthermore, we have the following strong convergences: 
$$ \int \rho_\epsilon{\vert u_\epsilon - \mathcal{R}(-t/\epsilon) w \vert^2}dx \rightarrow 0$$

and 
$$\sqrt{\frac{de^{V_\epsilon}}{\int d e^{V_\epsilon} dx}} \rightarrow \sqrt{\rho}$$
in $L^\infty_t(L^2_x) $.

\end{thm}

\begin{rque}
\begin{enumerate}
\item Instead of a cut-off of order $\frac{1}{\epsilon}$, we could have chosen any function $\xi(\epsilon)$ such that for some $q\in ]3/2,s-1[$:
\begin{equation*}
\begin{split}
\frac{1}{\xi(\epsilon)} \rightarrow_{\epsilon \rightarrow 0}& 0, \\
\xi(\epsilon)^{s-q-2} \epsilon  \rightarrow_{\epsilon \rightarrow 0}& 0.
 \end{split}
 \end{equation*}
 
 The choice $\xi(\epsilon)=\frac{1}{\epsilon}$ yields a sharp convergence rate.
 
 \item Concerning the rate of convergence, the proof below actually shows that for any $q\in ]3/2,s-1[$, if:
  \[
 \mathcal{H}_\epsilon(0) \leq C\epsilon^{\min (\alpha, s-q-1,1)},
 \]
 then there exists $C_q$ depending on $q$ such that locally uniformly in time:
 \[
 \mathcal{H}_\epsilon(t) \leq C_q \epsilon^{\min (\alpha, s-q-1,1) }.
 \]
 
\item When $s> 7/2$, we can observe in the proof that we actually do not need any cut-off in frequency. In this case the convergence is of order $\epsilon^{\min(1,\alpha)}$.

 \end{enumerate}
 
\end{rque}

\begin{proof}

We assume that $s\in ]5/2, 7/2]$. When $s>7/2$, the proof is actually much simpler, as we do not need any cut-off in frequency and all the estimates are straightforward.

\par
\textbf{Step 1} We first show that $\Vert z_\epsilon \Vert_{L^\infty_t([0,T], H^{s-1})} \leq  C$.  Let us observe that we do not use the cut-off in frequency here.
We have:
$$\int_{\mathbb{R}^3} \left( 1 + \vert \xi \vert^2 \right)^{s-1} \vert \mathcal{F}z^\rho_\epsilon(\xi)  \vert ^2 d\xi = \Vert z_\epsilon^\rho \Vert_{H^{s-1}}^2.$$

We then estimate:
\begin{eqnarray*}
& &\int_{\vert \xi \vert \leq\frac{1}{\epsilon}} \left( 1 + \vert \xi \vert^2 \right)^{s-1}\left\vert \mathcal{F}\left(\nabla_{x_\perp}. (\mathcal{S}(t/\epsilon)w) \right)\right\vert^2 d\xi  \\
&\leq& C \int_{\vert \xi \vert \leq\frac{1}{\epsilon}} \left( 1 + \vert \xi \vert^2 \right)^{s-1+1}  \left\vert \mathcal{F}\mathcal{S}(t/\epsilon) w\right\vert^2 d\xi \\
&\leq& C  \int_{\vert \xi \vert \leq\frac{1}{\epsilon}} \left( 1 + \vert \xi \vert^2 \right)^{s}  \left\vert \mathcal{F}\mathcal{S}(t/\epsilon)w\right\vert^2 d\xi \\
&\leq& C  \Vert \mathcal{S}(t/\epsilon) w \Vert_{H^{s}}^2 \leq C  \Vert  w \Vert_{H^{s}}^2  .
\end{eqnarray*}

Similarly we have:
\begin{equation}
\begin{split}
\int_{\vert \xi \vert \leq\frac{1}{\epsilon}} \left( 1 + \vert \xi \vert^2 \right)^{s-1} \vert \mathcal{F}\left(\nabla_x H . \mathcal{S}(t/\epsilon) w_\perp \right)\vert^2 d\xi \leq C \Vert w_\perp \Vert_{H^{s-1}}^2. 
\end{split}
\end{equation}

Finally we compute:
\begin{equation}
\begin{split}
&\int_{\mathbb{R}^3}\int_{\mathbb{R}^3}\left( 1 + \vert \xi \vert^2 \right)^{s-1} \mathbbm{1}_{\vert \xi- \eta \vert + \vert \eta \vert \leq\frac{1}{\epsilon}} \vert \mathcal{F}(\mathcal{S}(t/\epsilon) w_\perp)(\eta).\mathcal{F} \left(\nabla_{x_\perp}  \log\left(\frac{\rho}{d}\right)\right)(\xi-\eta)\vert^2 d\eta \\
\leq& C  \Vert \mathcal{S}(t/\epsilon) w  \nabla_x \log \frac{\rho}{d} \Vert_{H^{s-1}}^2  \leq C  \Vert  w \Vert_{H^{s-1}}^2 \Vert  \log\frac{\rho}{d}\Vert_{H^{s}}^2,
\end{split}
\end{equation}
since $H^{s-1}(\mathbb{R}^3)$ is an algebra. This proves that there exists a constant depending on $w$ and $\log \frac\rho d$:
$$\Vert z_\epsilon^{\rho} \Vert_{L^\infty_t([0,T], H^{s-1})} \leq  C. $$

Likewise, we prove that:
$$\Vert z_\epsilon^{w} \Vert_{L^\infty_t([0,T], H^{s-1})} \leq  C. $$

This yields that $\epsilon z_\epsilon \rightarrow 0$ in $L^\infty_t([0,T], H^{s-1})$. 

Now, let $q \in ]3/2,s-1[$ be a fixed parameter. We are interested in the $H^{q+1}$ norm of $z_\epsilon$.  Since $q+2 >7/2 \geq s$, the $H^{q+2}$ norm of $\log \frac \rho d$ and $w$ is not necessarily well-defined\footnote{When $s>7/2$, $q$ can be chosen such that $q+2\leq s$ and thus we have $\Vert z_\epsilon \Vert_{L^\infty_t([0,T], H^{q+1})}  \leq C$.} and we use this time the cut-off in frequency to lower down the regularity to $H^s$:
\begin{eqnarray*}
& &\int_{\vert \xi \vert \leq\frac{1}{\epsilon}} \left( 1 + \vert \xi \vert^2 \right)^{q+1}\left\vert \mathcal{F}\left(\nabla_{x_\perp}. (\mathcal{S}(t/\epsilon)w) \right)\right\vert^2 d\xi  \\
&\leq& C \int_{\vert \xi \vert \leq\frac{1}{\epsilon}} \left( 1 + \vert \xi \vert^2 \right)^{s}\left( 1 + \vert \xi \vert^2 \right)^{q+2-s}  \left\vert \mathcal{F}\mathcal{S}(t/\epsilon) w\right\vert^2 d\xi \\
&\leq& C/\epsilon^{q+2-s}  \int_{\vert \xi \vert \leq\frac{1}{\epsilon}} \left( 1 + \vert \xi \vert^2 \right)^{s}  \left\vert \mathcal{F}\mathcal{S}(t/\epsilon)w\right\vert^2 d\xi \\
&\leq&  C/\epsilon^{q+2-s} \Vert  w \Vert_{H^{s}}^2  .
\end{eqnarray*}

Treating the other terms similarly, we obtain:
\begin{equation}
\label{Hq+1-norm}
\Vert z_\epsilon \Vert_{L^\infty_t([0,T], H^{q+1})}  \leq C/\epsilon^{q+2-s}.
\end{equation}

\textbf{Step 2} We denote $X_\epsilon = \begin{pmatrix}{\rho}e^{\epsilon z^\rho_\epsilon} \\  w + \epsilon z^w_\epsilon \end{pmatrix}$. 
We introduce the filtered acceleration operator defined by:

\begin{equation}
B_\epsilon({\bar\rho}, \bar u)= \begin{pmatrix}\partial_t \log\left(\frac{\bar\rho}{d}\right) + \nabla_x. \mathcal{R}(-t/\epsilon) \bar u + \mathcal{R}(-t/\epsilon)\bar u. \nabla_x  \log\left(\frac{\bar\rho}{d}\right) - \nabla_x H . \mathcal{R}(-t/\epsilon)\bar u\\  \partial_t \bar u + \mathcal{R}(-t/\epsilon)\bar u. \nabla_x \bar u + \mathcal{R}(-t/\epsilon)\nabla_x \log\left(\frac{\bar\rho}{d}\right) \end{pmatrix}.
\end{equation}

We show that $X_\epsilon$ is an approximate zero of the filtered acceleration operator $B_\epsilon$ in the sense that (we recall that $q \in ]3/2,s-1[$):
$$\Vert B_\epsilon(X_\epsilon) \Vert_{L^\infty_t([0,T], H^q)} \rightarrow 0,$$ 
when $\epsilon \rightarrow 0$. 

By definition of $\mathcal{S}$, we have:
\[
\frac{d\mathcal{S}(t/\epsilon)}{dt}= \frac 1 \epsilon \mathcal{R}(-t/\epsilon)_\perp.
\]

Hence we have:
\begin{eqnarray*}
& &\partial_t \mathcal{F} \left(\log(\rho/d)+ \epsilon z_\epsilon^\rho\right) \\
&=& - \mathcal{F} \left(\partial_{x_\parallel} w_\parallel + w_\parallel \partial_{x_\parallel} \log(\rho/d) - \partial_{x_\parallel} H w_\parallel  \right) \\
&-& \mathbbm{1}_{\vert \xi \vert \leq \frac{1}{\epsilon}} \mathcal{F}\left(\nabla_{x_\perp}. (\mathcal{R}(-t/\epsilon)w)  - \nabla_x H . \mathcal{R}(-t/\epsilon) w_\perp \right)\\
&-&\int_{\mathbb{R}^3} \mathbbm{1}_{\vert \xi- \eta \vert + \vert \eta \vert \leq\frac{1}{\epsilon}} \mathcal{F}(\mathcal{R}(-t/\epsilon) w_\perp(\eta)) \mathcal{F} (\nabla_{x_\perp}  \log\left(\frac{\rho}{d}\right))(\xi-\eta)d\eta \\
&+& D_\epsilon(t,\xi),
\end{eqnarray*}
where $D_\epsilon$ is defined by:
\begin{eqnarray*}
D_\epsilon(t,\xi):= &-& \epsilon \mathbbm{1}_{\vert \xi \vert \leq \log(\epsilon)} \mathcal{F}\left(\nabla_{x_\perp}. (\mathcal{S}(t/\epsilon)\partial_t w)  - \nabla_x H . \mathcal{S}(t/\epsilon) \partial_t w_\perp \right) \\
&-& \epsilon \int_{\mathbb{R}^3} \mathbbm{1}_{\vert \xi- \eta \vert + \vert \eta \vert \leq\frac{1}{\epsilon}} \mathcal{F}\mathcal{S}(t/\epsilon) \partial_t w_\perp(\eta) \mathcal{F} \left(\nabla_{x_\perp}  \log\left(\frac{\rho}{d}\right)\right)(\xi-\eta)d\eta\\ 
&-&\epsilon \int_{\mathbb{R}^3} \mathbbm{1}_{\vert \xi- \eta \vert + \vert \eta \vert \leq\frac{1}{\epsilon}} \mathcal{F}\mathcal{S}(t/\epsilon) w_\perp(\eta) \mathcal{F} \left(\nabla_{x_\perp}  \partial_t \log\left(\frac{\rho}{d}\right)\right)(\xi-\eta)d\eta.
\end{eqnarray*}

Consequently, denoting by  $ B_{1,\epsilon}$ the operator in the first line of $B_\epsilon$ (resp. $B_{2,\epsilon}$ the second operator) we have:


\begin{eqnarray*}
\mathcal{F}  B_{1,\epsilon}(X_\epsilon)  =  T_{1,\epsilon}(t,\xi) + T_{2,\epsilon}(t,\xi) + D_\epsilon(t,\xi),
\end{eqnarray*}
with:
\begin{equation}
\begin{split}
 T_{1,\epsilon}(t,\xi)=& \ \mathbbm{1}_{\vert \xi \vert >  \frac{1}{\epsilon}} \mathcal{F}\left(\nabla_{x_\perp}. (\mathcal{R}(-t/\epsilon)w_\perp)  - \nabla_x H . \mathcal{R}(-t/\epsilon) w_\perp \right)\\
+&\int_{\mathbb{R}^3} \mathbbm{1}_{\vert \xi- \eta \vert + \vert \eta \vert >\frac{1}{\epsilon}} \mathcal{F}(\mathcal{R}(-t/\epsilon) w_\perp)(\eta) \mathcal{F} \left(\nabla_{x_\perp}  \log\left(\frac{\rho}{d}\right)\right)(\xi-\eta)d\eta,
 \end{split}
\end{equation}
\begin{equation}
\begin{split}
 T_{2,\epsilon}(t,\xi)=& \ \epsilon \mathcal{F}\left(\nabla_x. \mathcal{R}(-t/\epsilon) z^w_\epsilon + \mathcal{R}(-t/\epsilon)w. \nabla_x z^\rho_\epsilon +  \mathcal{R}(-t/\epsilon) z_\epsilon^w. \nabla_x \log\frac{\rho}{d} - \nabla_x H . \mathcal{R}(-t/\epsilon) z^w_\epsilon\right) \\
+& \epsilon^2 \mathcal{F}(\mathcal{R}(-t/\epsilon) z_\epsilon^w. \nabla_x  z_\epsilon^\rho).
 \end{split}
\end{equation}

\begin{rque}
Without corrector ($z_\epsilon=0$) we have:
 \begin{eqnarray*}
 T_{1,\epsilon}(t,\xi) &=&  \mathcal{F}\left(\nabla_{x_\perp}. (\mathcal{R}(-t/\epsilon)w_\perp)  + \nabla_x H . \mathcal{R}(-t/\epsilon) w_\perp \right)\\
&+&\int_{\mathbb{R}^3}  \mathcal{F}(\mathcal{R}(-t/\epsilon) w_\perp)(\eta) \mathcal{F} \left(\nabla_{x_\perp}  \log\left(\frac{\rho}{d}\right)\right)(\xi-\eta)d\eta 
\end{eqnarray*}
These terms only weakly but not strongly converge to $0$ as $\epsilon$ goes to $0$ : this is why we have to add the corrector.

When $z_\epsilon$ is defined without without cut-off in frequency, we notice that the have $T_{1,\epsilon}(t,\xi)=0$. 

\end{rque}

\textbf{Estimating $T_{1,\epsilon}$} 

We need the $H^s$ regularity of $w$ and $\log \rho/d$ in order to get some decay in $\epsilon$ for $T_{1,\epsilon}$, by using, for any $\beta>0$:
\begin{eqnarray*}
& &\mathbbm{1}_{\vert \xi \vert >  \frac{1}{\epsilon}}  \leq{(1+\vert \xi \vert^2)^{\beta}}{{\epsilon}^{2 \beta}},\\
& &  \mathbbm{1}_{\vert \xi- \eta \vert + \vert \eta \vert >\frac{1}{\epsilon}}  \leq 2{(\vert \xi -\eta\vert^2+\vert \eta \vert^2)^{\beta}}{ \epsilon^{2 \beta}}.
\end{eqnarray*}

Therefore we have:
\begin{equation*}
\begin{split}
 &\int_{\mathbb{R}^3} (1+\vert \xi \vert^2)^q\mathbbm{1}_{\vert \xi \vert >   \frac{1}{\epsilon}} \vert\mathcal{F}\left(\nabla_{x_\perp}. (\mathcal{R}(-t/\epsilon)w_\perp)\right)\vert^2d\xi\\
\leq& {C}{ \epsilon^{2s-2(q+1)}}\int_{\mathbb{R}^3} (1+\vert \xi \vert^2)^{q+1+s-(q+1)}  \vert\mathcal{F}\left( \mathcal{R}(-t/\epsilon)w_\perp \right)\vert^2d\xi\\
\leq& {C}{\epsilon^{2s-2(q+1)}} \Vert w \Vert_{H^s}^2.
\end{split}
\end{equation*}

We handle the other terms by the same method.There exists a constant $C$ depending only on $q$ and the $H^s$ norm of $\log \frac \rho d$ and $w$ such that:
\begin{equation}
\left(\int_{\mathbb{R}^3} (1+\vert \xi \vert^2)^q \vert T_{1,\epsilon} \vert^2 d\xi \right)^{1/2}\leq C \epsilon^{s-q-1}.
\end{equation}

\textbf{Estimating $T_{2,\epsilon}$} 

We can use estimate (\ref{Hq+1-norm}). As a result, there exists a constant $C$ depending only on $q$ and the $H^s$ norm of $\log \frac \rho d$ and $w$ such that:
\begin{equation}
\left(\int_{\mathbb{R}^3} (1+\vert \xi \vert^2)^q \vert T_{2,\epsilon} \vert^2 d\xi\right)^{1/2} \leq C( \epsilon \times \epsilon^{s-q-2} +\epsilon^2 \times \epsilon^{s-q-2} )\leq C \epsilon^{s-q-1}.
\end{equation}

\textbf{Estimating $D_\epsilon$} 

We only have $\partial_t w \in H^{s-1}$ and consequently we do not necessarily have $\partial_t w \in H^{q+1}$. Nevertheless, we can use the cut-off in frequency to lower the regularity down to only $H^{s-1}$:
\begin{equation*}
\begin{split}
 &\int_{\mathbb{R}^3} (1+\vert \xi \vert^2)^q\mathbbm{1}_{\vert \xi \vert \leq \frac{1}{\epsilon}} \vert \mathcal{F}\left(\nabla_{x_\perp}. (\mathcal{S}(t/\epsilon)\partial_t w) \right)\vert^2 d\xi \\
 \leq& C \int_{\mathbb{R}^3} (1+\vert \xi \vert^2)^{q+1}\mathbbm{1}_{\vert \xi \vert \leq\frac{1}{\epsilon} } \vert \mathcal{F}\left( (\mathcal{S}(t/\epsilon)\partial_t w) \right)\vert^2 d\xi \\
 \leq& C{\frac{1}{\epsilon^{2(q+1)-2(s-1)}}}\int_{\mathbb{R}^3}\mathbbm{1}_{\vert \xi \vert \leq \frac{1}{\epsilon}}  (1+\vert \xi \vert^2)^{s-1}\vert \mathcal{F}\left( (\mathcal{S}(t/\epsilon)\partial_t w) \right)\vert^2 d\xi \\
 \leq& C {\frac{1}{\epsilon^{2(q+1)-2(s-1)}}} \Vert \partial_t w \Vert_{H^{s-1}}.
\end{split}
\end{equation*}

Following the same method for the other terms, we finally obtain:

\begin{equation}
\left(\int_{\mathbb{R}^3} (1+\vert \xi \vert^2)^q \vert D_{\epsilon} \vert^2 d\xi\right)^{1/2} \leq C \epsilon^{s-1-q}.
\end{equation}

Gathering the pieces together, there exists a constant $C>0$ depending on $q$, the $H^s$ norm of $\log \frac \rho d$, $w$ and the $H^{s-1}$ norm of $\partial_t \log \frac \rho d$, $\partial_t w$ such that
\begin{eqnarray*}
\Vert  B_{1,\epsilon}(X_\epsilon) \Vert_{H^q} \leq C {\epsilon^{s-q-1}}.
\end{eqnarray*}

As a consequence, we have proved:

\begin{equation}
\Vert B_{1,\epsilon}(X_\epsilon) \Vert_{L^\infty_t([0,T], H^q)} \rightarrow 0.
\end{equation}

Arguing similarly for  $B_{2,\epsilon}$, we finally deduce that

$$\Vert B_\epsilon(X_\epsilon) \Vert_{L^\infty_t([0,T], H^q)} \rightarrow 0.$$


\textbf{Step 3} 
Finally we check that uniformly in time  :
$$ \int_0^t \int A_\epsilon(\bar\rho_\epsilon, \bar u_\epsilon). \begin{pmatrix} -m_\epsilon +\bar{\rho}_\epsilon \\ J_\epsilon -\rho_\epsilon\bar{u}_\epsilon\end{pmatrix} dx ds
\rightarrow 0,$$
as $\epsilon$ goes to $0$. We recall that $\bar \rho_\epsilon$ and $\bar u_\epsilon$ were defined in (\ref{defbar}).

First we have to check that
$$\Vert A_\epsilon(\bar\rho_\epsilon ,  \bar u_\epsilon) \Vert_{L^\infty_t([0,T], H^q)} \rightarrow 0.$$ 

This is clear in view of Step 2, since we have:
$$A_\epsilon(\bar\rho_\epsilon ,  \bar u_\epsilon) =\begin{pmatrix} B_{1,\epsilon}(X_\epsilon) \\ \mathcal{R}(t/\epsilon) B_{2,\epsilon}(X_\epsilon) \end{pmatrix}.$$ 
and $\mathcal{R}(t/\epsilon)$ is an isometry on any $H^s(\mathbb{R}^3)$.

We denote $A_\epsilon= \begin{pmatrix}A_{1,\epsilon} \\ A_{2,\epsilon}\end{pmatrix}$ and evaluate:
\begin{eqnarray*}
\left\vert\int_0^t \int A_{1,\epsilon}(\bar\rho_\epsilon, \bar u_\epsilon). \begin{pmatrix} -m_\epsilon +\bar{\rho}_\epsilon \end{pmatrix} dx ds\right\vert \leq \int_0^t \int \left\vert A_{1,\epsilon}(\bar\rho_\epsilon, \bar u_\epsilon)m_\epsilon  \right\vert dxds + \leq \int_0^t \int \left\vert A_{1,\epsilon}(\bar\rho_\epsilon, \bar u_\epsilon)\bar\rho_\epsilon \right\vert dxds \\
\leq C \Vert A_\epsilon(\bar\rho_\epsilon ,  \bar u_\epsilon) \Vert_{L^\infty_t([0,T], L^\infty)}\left(\Vert m_\epsilon \Vert_{L^\infty_t(L^1_x)}
+ \Vert \rho \Vert_{L^\infty_t(L^1_x)}\Vert e^{\epsilon z_\epsilon^\rho} \Vert_{L^\infty_{t,x}}\right) \\
\leq C \Vert A_\epsilon(\bar\rho_\epsilon ,  \bar u_\epsilon) \Vert_{L^\infty_t([0,T], H^q)}\left(1
+\Vert e^{\epsilon z_\epsilon^\rho} \Vert_{L^\infty_{t,x}}\right),
 \end{eqnarray*} 
 by Sobolev embedding, $H^q(\mathbb{R}^3) \rightarrow L^\infty(\mathbb{R}^3)$ ($q>3/2$). By the estimates of Step 1, there exists $C>0$ independent of $\epsilon$ such that:
 $$\Vert e^{\epsilon z_\epsilon^w} \Vert_{L^\infty_{t,x}}\leq C.$$
 
 In the other hand,
 
 \begin{eqnarray*}
& &\left\vert \int_0^t \int A_{2,\epsilon}(\bar\rho_\epsilon, \bar u_\epsilon). \begin{pmatrix}  J_\epsilon -\rho_\epsilon\bar{u}_\epsilon\end{pmatrix} dx ds\right\vert \\
 &\leq& C \Vert A_\epsilon(\bar\rho_\epsilon ,  \bar u_\epsilon) \Vert_{L^\infty_t([0,T], L^\infty)}\left(\Vert J_\epsilon \Vert_{L^\infty_t(L^1_x)}
+ \Vert \rho_\epsilon \Vert_{L^\infty_t(L^1_x)}\Vert \mathcal{R}(-t/\epsilon)(w + \epsilon z_\epsilon^w) \Vert_{L^\infty_{t,x}}\right)\\
&\leq& \Vert A_\epsilon(\bar\rho_\epsilon ,  \bar u_\epsilon) \Vert_{L^\infty_t([0,T], H^q)}\left(1
+\Vert w\Vert_{L^\infty_t(H^s)}+ \Vert \epsilon z_\epsilon^w \Vert_{L^\infty_{t,x}}\right)
 \end{eqnarray*}
and the conclusions follows.

One can also readily check, using (\ref{esti-G}), that $G_\epsilon(t) \rightarrow 0$ uniformly in time.

Finally this proves that $\mathcal{H}_\epsilon(t) \rightarrow 0$ uniformly in time, as soon as  $\mathcal{H}_\epsilon(0) \rightarrow 0$. 

Using the estimates of Step 1, we check that
$$\left(\rho e^{\epsilon z_\epsilon^\rho}, w+\epsilon z_\epsilon^w\right) \rightarrow (\rho, w)$$
in $L^\infty([0,T],L^\infty)$.

In order to apply Gronwall's inequality to the inequality (\ref{stability-B}), there remains to check that $\Vert \nabla_x \bar u_\epsilon \Vert_{L^\infty}$ is uniformly bounded in $\epsilon$. It is sufficient to check that $\Vert\epsilon \nabla_x  z_\epsilon^w \Vert_{L^\infty}$ is uniformly bounded. According to (\ref{Hq+1-norm}) and by Sobolev embedding $  H^{q} \rightarrow  L^\infty$, we have:
\begin{equation}
\Vert\epsilon \nabla_x  z_\epsilon^w \Vert_{L^\infty_{t,x}} \leq  \Vert\epsilon \nabla_x  z_\epsilon^w \Vert_{L^\infty_t H_x^{q}} \leq C \epsilon^{s-q-1} \leq C.
\end{equation}

Then, the other conclusions easily follow as in the end of the proof of Theorem \ref{theoS}.
\end{proof}

\section*{Annex}
\subsection*{Scaling of the Vlasov-Poisson systems ({S}), ({S'}) and ({L})}

Let us introduce the dimensionless variables and unknowns:

$$\tilde{t}= \frac t \tau \quad \tilde{x}= \frac x L \quad \tilde{v}= \frac v {v_{th}},$$
$$f(t,x,v)= \bar f \tilde{f}(\tilde{t},\tilde{x},\tilde{v}) \quad V(t,x)= \bar V \tilde V(\tilde{t},\tilde{x}) \quad E(t,x)=\bar E \tilde E(\tilde t, \tilde x). $$

 Then the Vlasov equation with Poisson equation (\ref{MB}) equation states:

   \begin{equation}
\left\{
    \begin{array}{ll}
  \partial_{\tilde t} \tilde f_\epsilon +  \frac{v_{th}\tau}{L}\tilde v.\nabla_{\tilde x}\tilde f_\epsilon + \frac{e\bar E \tau}{mv_{th}}\tilde E_\epsilon.\nabla_{\tilde v} \tilde f_\epsilon =0  \\
 \frac{\bar E L}{\bar V} \tilde E_\epsilon = -\nabla_{\tilde x} \tilde V_\epsilon\\
-\frac{\epsilon_0\bar V}{L^2}\Delta_{\tilde x} \tilde V_\epsilon =e\bar f v_{th}^3\int \tilde f_\epsilon d \tilde v - e\bar d  \tilde d e^{\frac{e\bar V}{k_B T_e}\tilde V_\epsilon}\\
     \tilde f_{\epsilon,\vert \tilde t=0} =\tilde f_{0,\epsilon}, \quad \bar f L^3 v_{th}^3\int \tilde f_{0,\epsilon} d \tilde vd \tilde x =1.\\
    \end{array}
  \right.
\end{equation}

In order to ensure that $\int \tilde f d\tilde{x} d\tilde{v}=1$, it is natural to set:
 $$\bar f  L^3 v_{th}^3= 1.$$

Moreover we consider the normalizations:

$$\frac{v_{th}\tau}{L}=1, \quad  \frac{\bar E L}{\bar V}=1, $$ 
$$\frac{e\bar V}{k_B T_e} =1, $$
$$\bar f v_{th}^3 = \bar d.$$

This implies that:
$$\frac{e\bar E \tau}{mv_{th}} =  \frac{v_{th} \tau}{L}=1.$$

Now we observe that:
$$\frac{\epsilon_0\bar V }{ e \bar f v_{th}^3 } =  \frac{\epsilon_0 k_B T_e }{e^2 \times 1/L^3} =  {\lambda_D^2}, $$
where $\lambda_D$ is the Debye length.

The quasineutral scaling consists in considering the ordering:
$$ \frac {\lambda_D^2} {L^2} =\epsilon,  $$
with $\epsilon$ a small parameter.

With this scaling we get the following dimensionless system of equations (we forget the $\tilde{}$ for the sake of readability):

   \begin{equation}
\left\{
    \begin{array}{ll}
  \partial_t f_\epsilon + v.\nabla_x f_\epsilon + E_\epsilon.\nabla_v f_\epsilon =0  \\
  E_\epsilon = -\nabla_x V_\epsilon\\
-\epsilon \Delta_x V_\epsilon = \int f_\epsilon dv - d e^{V_\epsilon}\\
    f_{\epsilon,\vert t=0} =f_{0,\epsilon} \geq 0, \int f_{0,\epsilon} dvdx =1.\\
    \end{array}
  \right.
\end{equation}
This system is nothing but System (S). We get Systems (S') and (L) with the same nondimensionalization.

\begin{rque}
To be rigorous we should also consider the confinement force $-\nabla_x H$ on the ions, but we will not do so for the sake of simplicity and readability. Nevertheless we could handle such an external force with only minor changes in the following.
\end{rque}

\subsection*{Scaling of the Vlasov-Poisson equation (\ref{(S')-B})}

We once again consider the nondimensionalization analysis of system (S'), this time including the magnetic field:
$$B = \bar B e_\parallel.$$

This yields:

   \begin{equation}
\left\{
    \begin{array}{ll}
  \partial_{\tilde t} \tilde f_\epsilon +  \frac{v_{th}\tau}{L}\tilde v.\nabla_{\tilde x}\tilde f_\epsilon +\left(\frac{e\bar E \tau}{mv_{th}}\tilde E_\epsilon+ \frac{e \bar B}{m} \tau  \tilde v \wedge e_\parallel\right).\nabla_{\tilde v} \tilde f_\epsilon =0  \\
 \frac{\bar E L}{\bar V} \tilde E_\epsilon = -\nabla_{\tilde x} \tilde V_\epsilon\\
-\frac{\epsilon_0\bar V}{L^2}\Delta_{\tilde x} \tilde V_\epsilon =e\bar f v_{th}^3\int \tilde f_\epsilon d \tilde v - e\frac{\tilde d e^{\frac{e\bar V}{k_B T_e}\tilde V_\epsilon}}{\int \tilde d e^{\frac{e\bar V}{k_B T_e}\tilde V_\epsilon}dx}\\
     \tilde f_{\epsilon,\vert \tilde t=0} =\tilde f_{0,\epsilon}, \quad \bar f L^3 v_{th}^3\int \tilde f_{0,\epsilon} d \tilde vd \tilde x =1.\\
  \end{array}
  \right.
\end{equation}

We set $\Omega=  \frac{e\bar B}{m}$ : this is the cyclotron frequency (also referred to as the gyrofrequency).
We also consider the so-called electron Larmor radius (or electron gyroradius) $r_L$ defined by:

\begin{equation}
r_L=\frac{ v_{th}}{\Omega}= \frac{m v_{th}}{e \bar B}.
\end{equation}
This quantity can be physically understood as the typical radius of the helix around axis $e_\parallel$ that the particles follow, due to the intense magnetic field.

The Vlasov equation now reads:
$$ \partial_{\tilde t} \tilde f_\epsilon +  \frac{r_L}{L}\Omega \tau\tilde v.\nabla_{\tilde x}\tilde f_\epsilon  +\left(\frac{\bar E }{\bar B v_{th}}\Omega \tau\tilde E_\epsilon+ \Omega \tau   \tilde v \wedge e_\parallel\right).\nabla_{\tilde v} \tilde f_\epsilon =0 .
$$

The "strong magnetic field" ordering consists in setting:
$$\Omega \tau =\frac 1 \epsilon, \quad \frac{\bar E }{\bar B v_{th}}=\epsilon,$$
$$\frac{r_L}{L}=\epsilon.$$

The quasineutral scaling we consider is:
$$ \frac {\lambda_D^2} {L^2} =\epsilon^{2\alpha} ,$$
with $\alpha>0$. From the physical point of view, it means that we consider that both the Larmor radius and the Debye length vanish. We observe that:

$$\frac {\lambda_D}{r_L}= \epsilon^{\alpha-1}$$

In most practical situations, $\lambda_D << r_L$, so that  the range of parameters $\alpha>1$ is particularly physically relevant. Finally, by having the same normalizations as before, we get in the end:

\begin{equation*}
\left\{
    \begin{array}{ll}
  \partial_t f_{\epsilon} + v.\nabla_x f_{\epsilon} + \left(E_{\epsilon}+ \frac{v\wedge b}{\epsilon}\right).\nabla_v f_{\epsilon} =0  \\
  E_{\epsilon} = -\nabla_x V_{\epsilon}\\
-\epsilon^{2\alpha} \Delta_x V_{\epsilon} = \int f_{\epsilon} dv -\frac{d e^{V_{\epsilon}}}{\int d e^{V_\epsilon}dx}
\\
    f_{{\epsilon},\vert t=0} =f_{0,\epsilon}, \quad \int f_{0,\epsilon} dvdx =1.\\
    \end{array}
  \right.
\end{equation*}

\bibliographystyle{plain}
\bibliography{quasineutral}

\begin{thebibliography}{10}

\bibitem{Ar}
A.A. {Arsenev}.
\newblock {Existence in the large of a weak solution of Vlasov's system of
  equations}.
\newblock {\em Z. Vychisl. Mat. Mat. Fiz}, 15:136--147, 1975.

\bibitem{BOS}
M.~Bostan.
\newblock {The Vlasov-Maxwell system with strong initial magnetic field:
  guiding-center approximation}.
\newblock {\em {Multiscale Model. Simul.}}, 6(3):1026--1058, 2007.

\bibitem{Bou}
F.~Bouchut.
\newblock {Global weak solution of the Vlasov-Poisson system for small
  electrons mass}.
\newblock {\em {Comm. Partial Differential Equations}}, 16:1337--1365, 1991.

\bibitem{BGP}
F.~{Bouchut}, F.~{Golse}, and M.~{Pulvirenti}.
\newblock {\em {Kinetic Equations and Asymptotic Theory}}.
\newblock Series in Applied Mathematics. Gauthier-Villars, 2000.

\bibitem{Br}
Y.~{Brenier}.
\newblock {Convergence of the Vlasov-Poisson system to the incompressible Euler
  equations}.
\newblock {\em Comm. Partial Differential Equations}, 25:737--754, 2000.

\bibitem{BL}
Y.~Brenier and G.~Loeper.
\newblock {A geometric approximation to the Euler equations: the
  Vlasov-Monge-Ampère system}.
\newblock {\em {Geom. Funct. Anal.}}, 14(6):1182--1218, 2004.

\bibitem{BMP}
Y.~Brenier, N.~Mauser, and M.~Puel.
\newblock {Incompressible Euler and e-MHD as scaling limits of the
  Vlasov-Maxwell system}.
\newblock {\em {Commun. Math. Sci.}}, 1(3):437--447, 2003.

\bibitem{CGPSR}
C.~Cheverry, I.~Gallagher, T.~Paul, and L.~{Saint-Raymond}.
\newblock {Trapping Rossby waves}.
\newblock {\em {C. R. Math. Acad. Sci. Paris}}, 347(15-16):879--884, 2009.

\bibitem{CG}
Y.~Cordier and E.~Grenier.
\newblock {Quasineutral limit of an Euler-Poisson system arising from plasma
  physics}.
\newblock {\em {Comm. Partial Differential Equations}}, 25(5):1099--1113, 2000.

\bibitem{DA}
R.L. Dewar and R.F. Abdullatif.
\newblock {Zonal flow generation by modulational instability}.
\newblock {\em {Proceedings of the CSIRO/COSNet Workshop on Turbulence and
  Coherent Structures, Canberra, Australia}}, World Scientific, eds. J.P.
  Denier and J.S. Frederiksen, 2006.

\bibitem{DPL}
R.J. Diperna and P.L. Lions.
\newblock {Solutions globales d'equations du type Vlasov-Poisson}.
\newblock {\em {C. R. Acad. Sci. Paris Ser. I Math.}}, 307(12):655--658, 1988.

\bibitem{Fitz}
R.~Fitzpatrick.
\newblock {The Physics of Plasmas}.
\newblock {\em {available at
  http://farside.ph.utexas.edu/teaching/plasma/380.pdf}}.

\bibitem{FS2}
E.~{Frénod} and E.~{Sonnendrücker}.
\newblock {The Finite Larmor Radius Approximation}.
\newblock {\em SIAM J. Math. Anal.}, 32(6):1227--1247, 2001.

\bibitem{GSR2}
F.~{Golse} and L.~{Saint-Raymond}.
\newblock {The Vlasov-Poisson system with strong magnetic field in quasineutral
  regime}.
\newblock {\em Mathematical Models and Methods in Applied Sciences},
  13(5):661--714, 2003.

\bibitem{GJV}
T.~Goudon, P.E. Jabin, and A.~Vasseur.
\newblock {Hydrodynamic limit for the Vlasov-Navier-Stokes equations. II. Fine
  particles regime.}
\newblock {\em {Indiana Univ. Math. J.}}, 53(6):1517--1536, 2004.

\bibitem{Gre5}
E.~Grenier.
\newblock {Oscillatory perturbations of the Navier-Stokes equations}.
\newblock {\em {J. Math. Pures Appl.}}, 76(9):477--498, 1997.

\bibitem{Gre3}
E.~{Grenier}.
\newblock {Pseudo-Differential Energy Estimates of Singular Perturbations}.
\newblock {\em Comm. Pure App. Math.}, 50:0821--0865, 1997.

\bibitem{DHK1}
D.~{Han-Kwan}.
\newblock {The three-dimensional finite Larmor radius approximation}.
\newblock {\em Asymptot. Anal.}, 66(1):9--33, 2010.

\bibitem{HM1}
A.~Hasegawa and K.~Mima.
\newblock {Pseudo-three-dimensional turbulence in magnetized nonuniform
  plasma}.
\newblock {\em {Phys. Fluids}}, 21(1):87--92, 1978.

\bibitem{HH}
E.~Horst and E.~Hunze.
\newblock {Weak solutions of the initial value problem for the unmodified
  nonlinear Vlasov equation}.
\newblock {\em {Math. Methods Appl. Sci.}}, 6:262--279, 1982.

\bibitem{LL}
M.A. Lieberman and A.J. Lichtenberg.
\newblock {Principles of plasma discharge and materials processing}.
\newblock {\em {New-York, Wiley}}, 1994.

\bibitem{LP}
P.L. Lions and B.~Perthame.
\newblock {Propagation of moments and regularity for the three-dimensional
  Vlasov-Poisson system}.
\newblock {\em {Invent. Math.}}, 105:415--430, 1991.

\bibitem{Loe}
G.~Loeper.
\newblock {Uniqueness of the solution to the Vlasov-Poisson system with bounded
  density}.
\newblock {\em {J. Math. Pures Appl.}}, 86:68--79, 2006.

\bibitem{Maj}
A.~Majda.
\newblock {Compressible fluid flow and systems of conservation laws in several
  space variables}.
\newblock {\em {Springer-Verlag}}, 1984.

\bibitem{Mas}
N.~{Masmoudi}.
\newblock {From Vlasov-Poisson system to the incompressible Euler system}.
\newblock {\em Comm. Partial Differential Equations}, 26(9):1913--1928, 2001.

\bibitem{MV}
C.~Mouhot and C.~Villani.
\newblock {On Landau damping}.
\newblock {\em Preprint}, 2009.

\bibitem{Pfa}
K.~Pfaffelmoser.
\newblock {Global classical solutions of the Vlasov-Poisson system in three
  dimensions for general initial data}.
\newblock {\em {J. Diff. Equations. }}, 95:281--303, 1992.

\bibitem{PSR}
M.~Puel and L.~Saint-Raymond.
\newblock {Quasineutral limit for the relativistic Vlasov-Maxwell system}.
\newblock {\em {Asymptot. Anal.}}, 40(3-4):303--352, 2004.

\bibitem{SR}
L.~Saint-Raymond.
\newblock {Hydrodynamic limits: some improvements of the relative entropy
  method}.
\newblock {\em { Ann. Inst. H. Poincare Anal. Non Lineaire}}, 26(3):705--744,
  2009.

\bibitem{Scha}
J.~Schaeffer.
\newblock {Global existence of smooth solutions to the Vlasov-Poisson system in
  three dimensions}.
\newblock {\em {Comm. Partial Differential Equations}}, 16(8-9):1313--1335,
  1991.

\bibitem{Sch}
S.~{Schochet}.
\newblock {Fast singular limits of hyperbolic PDEs}.
\newblock {\em J. Differential Equations}, 114:476--512, 1994.

\end{thebibliography}

\end{document}